%% file: Asterisque-intro-2.tex
\documentclass[11pt]{amsart}
\usepackage{amsmath,amssymb,amsfonts}
\usepackage{enumerate}
\usepackage{amscd, amssymb, latexsym, amsmath, amscd}
\usepackage{mathtools}
\usepackage[all]{xy}
\usepackage{tikz}
\usepackage{tikz-cd}
\usepackage{pb-diagram}
\usepackage{picture}
\usepackage{hyperref}
\usepackage{multicol,lipsum}

\newtheorem{thm}[equation]{Theorem}
\newtheorem{prop}[equation]{Proposition}
\newtheorem{cor}[equation]{Corollary}

\newtheorem{lemma}[equation]{Lemma}

\numberwithin{equation}{section}

\include{MathIncludes}

\begin{document}

\title[ L-groups and the Langlands program for covering groups]{L-groups and the Langlands program \\
 for covering groups:  \\
 a historical introduction}
 \author{Wee Teck Gan, Fan Gao and Martin H. Weissman}

\address{Wee Teck Gan: Department of Mathematics, National University of Singapore, 10 Lower Kent Ridge Road
Singapore 119076} 
\email{matgwt@nus.edu.sg}

\address{Fan Gao: Department of Mathematics, Purdue University, 150 N. University Street, West Lafayette, IN 47907}
\email{fan166@purdue.edu}

\address{Martin Weissman: Department of Mathematics, University of California, Santa Cruz, CA 95064, and Yale-NUS College, 6 College Ave East,  \# B1-01, Singapore 138614}
\email{weissman@ucsc.edu}

\maketitle

In this joint introduction to the present Asterisque volume, we shall give a short discussion of the historical developments in the study of nonlinear covering groups, touching on their structure theory, representation theory and the theory of automorphic forms. This serves as a historical motivation and sets the scene for the papers in this volume. Our discussion is necessarily subjective and will undoubtedly leave out the contributions of many authors, to whom we apologize in earnest. 
 
 \section{ Generalities}

A {\bf locally compact group} will mean a locally compact, Hausdorff, second countable topological group.  Let $G$ be a locally compact group and $A$ a locally compact abelian group.  We are interested in central extensions of $G$ by $A$. Let us first define this notion; our treatment in this section follows the classic paper of Moore \cite{Mo}.

\subsection{ Definition:} A {\bf central extension of $G$ by $A$} is a short exact sequence:
\[ \begin{CD}
1 @>>> A @>i>> E @>p>> G @>>> 1 \end{CD} \]
such that
\begin{itemize}
\item $E$ is a locally compact group;

\item $i$ is continuous and $i(A)$ is a closed subgroup of the center of $E$;

\item $p$ is continuous and induces a topological isomorphism $E/i(A) \isom G$. 
\end{itemize}

\noindent Equivalently, the third condition above can be replaced by the 
requirement 
that $p$ is continuous and open (cf. \cite[p.96]{Mi}). We will ultimately be interested in the case when $A$ is finite.

\subsection{ Definition:}  Let $E_1$ and $E_2$ be two extensions of $G$ by $A$.  An {\bf equivalence} from $E_1$ to $E_2$ is a continuous homomorphism $\phi \From E_1 \To E_2$ inducing the identity maps on $A$ and $G$:  
\[ \begin{CD}
1 @>>> A @>>> E_1 @>{p_1}>> G @>>> 1  \\
  @.   @V{=}VV     @V{\phi}VV    @V{=}VV @.   \\
1 @>>> A  @>>> E_2 @>{p_2}>> G @>>> 1 . \end{CD} \]
By the open mapping theorem, an equivalence is necessarily a topological isomorphism.  

Let the set of equivalence classes of central extensions of $G$ by $A$ be denoted by $\CExt(G,A)$. The set $\CExt(G,A)$ has a natural 
abelian group structure, as we now explain. 

Given two extensions $E_1$ and $E_2$ of $G$ by $A$, we set
\[  E =  \{ (h_1, h_2) \in E_1 \times E_2: p_1(h_1) = p_2(h_2) \}/ \delta(A) \]
where $\delta(a) = (a, a^{-1})$ is the skew  diagonal embedding. This is the quotient of the fiber product $E_1 \times_G E_2$ by the skew diagonal 
embedding.  Then $E$ is a central extension of $G$ by $A$,
\[ \begin{CD}
1 @>>> A @>i>>  E @>p>>G @>>> 1, \end{CD} \]  
by defining $i(a) = (a, 1) = (1,a) \in E$ and $p(h_1, h_2) = p_1(h_1) = p_2(h_2)$. 

This $E$ is the so-called Baer sum of $E_1$ and $E_2$, written $E_1 \Baer E_2$ and this operation makes $\CExt(G,A)$ into an abelian group.  In other words, the equivalence class of $E$ depends only on the equivalence classes of $E_1$ and $E_2$.  

In the context of abstract groups, the abelian group $\CExt(G,A)$ was first studied by Schur (1904) who introduced the notion of Schur multipliers. In modern language, Schur had introduced the cohomology group $H^2(G,A)$.  We will however not go so far back in time in our historical discussion; a modern survey of the central extensions of finite groups of Lie type can be found in \cite{Pr}.

\subsection{ Categorical point of view.}

If we fix $G$ and $A$ as before, define $\Cat{CExt}(G,A)$ to be the category whose objects are central extensions of $G$ by $A$, and whose morphisms are equivalences.  Since all equivalences are isomorphisms, the category $\Cat{CExt}(G,A)$ is a groupoid.  The Baer sum is functorial,
$$\Baer \From \Cat{CExt}(G,A) \times \Cat{CExt}(G,A) \To \Cat{CExt}(G,A),$$
making the category $\Cat{CExt}(G,A)$ into a (strictly commutative) Picard category \cite[D\'efinition 1.4.2]{D3}

The neutral object in this category is the direct product $G \times A$.  Given an object $E \in \Cat{CExt}(G,A)$, a {\bf splitting} of $E$ is an equivalence (i.e., a morphism) from $E$ to $G \times A$.  

If $j \From H \To G$ is a continuous homomorphism of locally compact groups, and $(E,i,p) \in \Cat{CExt}(G,A)$, then we may {\bf pull back} the extension $E$ to define
$$j^\ast E = \{ (h,e) \in H \times E : j(h) = p(e) \}.$$
Then $j^\ast E \in \Cat{CExt}(H,A)$ by defining $i' \From A \To j^\ast E$ by $i'(a) = (1,i(a))$ and $p'(h,e) = p(e)$.

If $f \From A \To B$ is a continuous homomorphism of locally compact abelian groups, we may {\bf push out} the extension $(E,i,p)$ to define
$$f_\ast E = (B \times E) / \overline{ \langle (f(a), i(a)^{-1} ) : a \in A \rangle }.$$
Typically, $f$ will be a closed map, and so it will not be necessary to take the closure in the quotient above.  Then $f_\ast E \in \Cat{CExt}(G,B)$ by defining $i'' \From B \To f_\ast E$ by $i''(b) = (b,1)$ and $p'' \From f_\ast E \To G$ by $p''(b,e) = p(e)$.

Pullback and pushout define additive functors of Picard categories,
$$f_\ast \From \Cat{CExt}(G,A) \To \Cat{CExt}(G,B), \quad f^\ast \From \Cat{CExt}(G,A) \To \Cat{CExt}(H,A).$$
For isomorphism classes, these define homomorphisms of abelian groups,
$$f_\ast \From \CExt(G,A) \To \CExt(G,B), \quad f^\ast \From \CExt(G,A) \To \CExt(H,A).$$


\subsection{ Cohomological interpretation.}
After the foundational work of Mackey \cite{Mac}, Moore wrote a series of papers \cite{Mo1-2, Mo3,Mo4} developing a cohomology theory 
for topological groups analogous to that for abstract groups. We summarize some of their results.

Moore defines for each $n \geq 0$ a cohomology group $H^n(G,A)$ using {\em measurable} cochains.  These groups are functors which are covariant in $A$ and contravariant in $G$.  Note however that since the category of locally compact abelian groups is not an abelian category, this cohomology theory is not a derived functor cohomology theory.  
We describe the low degree cohomology groups concretely.  Note that we are only interested in the case where $A$ is trivial as a $G$-module.  The 0-th cohomology group is $H^0(G,A) = A$.  The first cohomology $H^1(G,A)$ is the group of continuous homomorphisms $G \To A$. 

We describe $H^2(G,A)$ in more detail.  Let $Z^2(G,A)$ be the group of measurable normalized 2-cocycles $z \From G \times G \To A$; this means that $z(g,1) = z(1,g) = 1$ for all $g \in G$, and 
$$z(g_1 g_2, g_3) z(g_1, g_2) = z(g_1, g_2 g_3) z(g_2, g_2) \text{ for all } g_1, g_2, g_3 \in G.$$
Let $C^1(G,A)$ be the group of normalized 1-cochains:  measurable functions from $G$ to $A$ such that $f(1) = 1$.   If $f \in C^1(G,A)$ is a measurable function, its coboundary $\partial f \in Z^2(G,A)$ is defined by
$$\partial f(g_1, g_2) = f(g_2) \cdot f(g_1 g_2)^{-1} \cdot  f(g_1)^{-1}.$$

The resulting cohomology group $H^2(G,A) = Z^2(G,A) / \partial C^1(G,A)$ is naturally isomorphic to $\CExt(G,A)$.  

This can be understood categorically as follows.  Consider the (small, strictly commutative Picard) category $\Cat{H}^2(G,A)$, with objects set $Z^2(G,A)$, and where a morphism $z_1 \To z_2$ is defined to be an element $c \in C^1(G,A)$ such that $z_2 = z_1 + \partial c$.  The Picard category structure arises from the abelian group structures on $Z^2(G,A)$ and $C^1(G,A)$.  The isomorphism classes in $\Cat{H}^2(G,A)$ form the cohomology group $H^2(G,A)$.

Describe a functor from $\Cat{H}^2(G,A)$ to $\Cat{CExt}(G,A)$ as follows:  for an object $z \in Z^2(G,A)$, define an extension of $G$ by $A$ by $E = G \times A$, with multiplication
\[ (g_1,a_1) \cdot (g_2, a_2) = (g_1g_2, a_1a_2 \cdot z(g_1, g_2)), \]
and maps
\[  i: a \mapsto (1,a)\in E  \quad \text{and} \quad p: (g,a) \mapsto g \in G. \]

A theorem of Mackey \cite[Th\'eor\`eme 2]{Mac} gives $E$ a natural topology such that the above defines a locally compact group, and an extension of $G$ by $A$.  If $c \From z_1 \To z_2$ is a morphism in $\Cat{H}^2(G,A)$, i.e., $z_2 = z_1 + \partial c$, then $c$ defines an equivalence of central extensions $E_1 \To E_2$ by the formula $f(g,a) = (g, c(g) \cdot a)$.  The work of Mackey and Moore implies that this gives an equivalence of Picard categories, which we like to call ``incarnation'':
$$\Cat{Inc} \From \Cat{H}^2(G, A) \To \Cat{CExt}(G,A).$$
A consequence is the isomorphism of abelian groups, $H^2(G,A) \isom \CExt(G,A)$.

Surjectivity of this isomorphism is obtained as follows.  Given a central extension $A \Into E \Onto G$, Mackey proves that one can find a measurable section $s \From  G \To E$ (i.e.~so that $p \circ s = id$).  This is the best one can hope for: one cannot find a continuous section in general.  From $s$, one defines a measurable 2-cocycle by:
\[  z(g_1, g_2) = s(g_1 g_2) \cdot s(g_1)^{-1} s(g_2)^{-1}. \]
The map $(g,a) \mapsto s(g) \cdot i(a)$ gives an isomorphism from $\Cat{Inc}(z)$ to $E$.


\subsection{ $\CExt(G,-)$ as a moduli functor.} 
For another perspective, fix a locally compact group $G$.  The assignment $A \mapsto \CExt(G,A)$ gives a functor,
$$\CExt(G, -) \From \Cat{LCA} \To \Cat{Ab},$$
where $\Cat{LCA}$ denotes the category of locally compact abelian groups and $\Cat{Ab}$ denotes the category of (abstract) abelian groups. Indeed, we have seen above that $\CExt(G,A)$ has a natural abelian group structure and functoriality comes from pushout.  Given $f \From A \To B$ in $\Cat{LCA}$, pushout defines a group homomorphism $f_\ast \From \CExt(G,A) \To \CExt(G,B)$.

%


Regarding $\CExt(G,-)$ as a functor $\Cat{LCA} \To \Cat{Ab}$, it is natural to ask:

\noindent{\bf Question:} Is the functor $\CExt(G,-)$ representable? If so, describe the representing object $\pi_1(G)$ of $\Cat{LCA}$ explicitly.

As we shall see in the next section, this is the central motivating question behind the initial study of central extensions, as developed by Steinberg \cite{S1}, Moore \cite{Mo}, Matsumoto \cite{Ma}, Raghunathan-Prasad \cite{PR1, PR2, PR3} and others.  

If $\pi_1(G)$ exists, we call it the {\bf fundamental group} of $G$, in which case we have isomorphisms, functorial in $A$:
\[  H^2(G,A) \isom \CExt(G,A) \isom \Hom(\pi_1(G), A).  \]
Note that this functorial isomorphism can only be unique up to automorphisms of $\pi_1(G)$.  Also, observe that if $A = S^1$ (the unit circle), then
\[  \pi_1(G) \cong \text{ the Pontryagin dual $H^2(G, S^1)^{\vee}$ of $H^2(G, S^1)$.} \]

\subsection{ Universal extensions.}
Suppose that $\pi_1(G)$ exists. Then since 
\[  \CExt(G, \pi_1(G)) \cong \Hom(\pi_1(G), \pi_1(G)), \]
there is an element $\tilde{G}$ of $\CExt(G,\pi_1(G))$ 
corresponding to the identity automorphism of $\pi_1(G)$. This extension
\[  \begin{CD} 
1 @>>> \pi_1(G) @>>> \tilde{G} @>>> G @>>> 1 \end{CD} \]
is called a {\bf universal central extension} because it has the following universal property: given any central extension $E$ of $G$ by $A$, there exists a {\em unique} continuous homomorphism $\phi \From \tilde G \To E$ lying over the identity on $G$:
\[ \begin{CD}
1 @>>> \pi_1(G) @>>> \tilde G @>>> G @>>> 1  \\
  @.   @V{f \vert_{\pi_1(G)}}VV     @V{f}VV    @V{=}VV @.   \\
1 @>>> A  @>>> E @>>> G @>>> 1 . \end{CD} \]
The existence of such a universal central extension is equivalent to the representability of $\CExt(G,-)$.

\subsection{ Condition for representability.} 
Now let us consider the question about existence of $\pi_1(G)$. There is an obvious necessary condition for this existence. Indeed,
consider the trivial extension $G \times S^1$ where $S^1$ is the unit circle. If $\pi_1(G)$ exists and $\tilde{G}$ is a universal central 
extension,  then we have a unique map of extensions
\[ f: \tilde{G} \longrightarrow G \times S^1. \]
However, if $\phi: G \longrightarrow S^1$ is any continuous homomorphism, then the map
\[  f_{\phi} : \tilde{g} \mapsto f(\tilde g) \cdot (1, \phi(p(\tilde{g})) \]
is another morphism of extensions. Thus the uniqueness of $f$ implies that $\Hom(G, S^1) = 0$. 
In other words, if $\pi_1(G)$ exists, then $H^1(G, S^1) = 0$, or 
equivalently, $[G,G]$ is dense in $G$, in which case we say that $G$ is 
{\bf topologically perfect}.


One may ask if the necessary condition above is sufficient for the existence of $\pi_1(G)$? Moore has given examples to show that
it is not in general. We highlight some positive results in this direction due to Moore  \cite{Mo}:

\begin{prop} \label{P:suff}  In the following cases, $\pi_1(G)$ exists:
\begin{enumerate}
\item[(i)]
$G$ is a discrete group which is perfect (equivalently, topologically perfect, since the topology is discrete);
\item[(ii)] $G$ is topologically perfect and $H^2(G, S^1)$ is finite;
\item[(iii)] the component group $G / G^\circ$ of $G$ is compact and $G = [G,G]$ (i.e., $G$ is perfect).
\end{enumerate}
\end{prop}
Recall that when $\pi_1(G)$ exists, $\pi_1(G) \cong H^2(G, S^1)^{\vee}$.


\subsection{ Relative fundamental groups.}
One may also  consider the problem of classifying the central extensions of $G$ which are split over a subgroup $H \subset G$. 
Here, $H$ need not be a closed subgroup; there are applications in which $H$ is even dense in $G$. 
Thus we are interested in the representability of the functor $\Cat{LCA} \To \Cat{Ab}$ given by
\[  A \mapsto \Ker(H^2(G,A) \longrightarrow H^2(H,A)). \]
One has the following result \cite[Lemma 2.8]{Mo}:

\begin{prop} \label{P:relative}
Suppose that $\pi_1(G)$ and $\pi_1(H)$ both exist. The map $i: H \longrightarrow G$ induces $i_\ast: \pi_1(H) \longrightarrow \pi_1(G)$.
Define
\[  \pi_1(G,H) = \pi_1(G)/ \overline{i_\ast(\pi_1(H))}. \]
Then there is an isomorphism, functorial in $A$:
\[  \Ker(H^2(G,A) \longrightarrow H^2(H,A)) \cong \Hom(\pi_1(G,H), A). \]
In other words, the above functor is represented by $\pi_1(G,H)$. 
\end{prop}

We call $\pi_1(G,H)$ the fundamental group of $G$ relative to $H$.

\subsection{ Restricted direct product.}

Given a countable collection $(G_v, K_v)$, $v \in S$, where each $G_v$ is a locally compact group and $K_v$ is an open compact subgroup of $G_v$, one may form the restricted direct product 
\[  G = \prod_v' (G_v, K_v) \]
which is still a locally compact topological group.

The following proposition (based on \cite[Theorem 12.1]{Mo}) gives natural conditions under which $\pi_1(G)$ exists.


\begin{prop} \label{P:product}
Write $i_v \From K_v \Into G_v$ for the inclusions.  Assume that
\begin{itemize}
\item $\pi_1(G_v)$ exists for all $v$;

\item $\pi_1(K_v)$ exists for almost all $v$; 

\item $i_{v \ast}(\pi_1(K_v))$ is open in $\pi_1(G_v)$ for almost all $v$.
\end{itemize}
Then $\pi_1(G)$ exists and is equal to the restricted direct product
\[  \pi_1(G) = \prod_v' \left( \pi_1(G_v), i_{v*}(\pi_1(K_v)) \right). \]
\end{prop}

\noindent This concludes our discussion on generalities about central extensions and fundamental groups.


\section{ Abstract Chevalley Groups}

In this section, we shall specialise to the case where $G = \alg{G}(k)$, with $\alg{G}$ a connected reductive group over a field $k$.  In particular, when $\alg{G}$ is split, we give a summary of the beautiful work of Steinberg \cite{S1, S2}, Moore \cite{Mo} and Matsumoto \cite{Ma} which describes the fundamental group of an abstract Chevalley group.  


To begin, let $\alg{G}$ be a split, simple and  simply connected linear algebraic group over an infinite field $k$.  Set $G = \alg{G}(k)$ regarded as an abstract group (with discrete topology).  In this case, $G$ is known to be perfect. Thus $\pi_1(G)$ exists by Proposition \ref{P:suff}(i) and there is a universal central extension
\[  \begin{CD}
1 @>>> \pi_1(G) @>>> E_G @>>> G @>>> 1. \end{CD} \]

\subsection{ Steinberg's construction of $E_G$.}
In \cite{S1}, Steinberg gave an explicit construction of the universal central extension $E_G$ using generators and relations. 
Fix 
\[  \alg{T} \subset \alg{B} = \alg{T} \cdot \alg{U}^+ \subset \alg{G}, \]
a split maximal torus contained in a Borel subgroup of $\alg{G}$. This gives rise to a root system $\Phi$ with a set $\Delta$ of simple roots.  For each $\alpha \in \Phi$, one has a root subgroup $\alg{U}_{\alpha} \cong \alg{G}_a$.  If an isomorphism $x_{\alpha} \From \alg{U}_{\alpha} \xrightarrow{\sim} \alg{G}_a$ is chosen for every $\alpha \in \Phi$, we define families of elements of $G$ by
\[  \begin{cases}
w_{\alpha}(t)  = x_{\alpha}(t) x_{-\alpha}(-t^{-1})x_{\alpha}(t), \\
h_{\alpha}(t) = w_{\alpha}(t) w_{\alpha}(1)^{-1} \end{cases} \]
for $t \in k^{\times}$.

In \cite[Chapter 6]{S2}, Steinberg demonstrates that there exist a family of isomorphisms $\{ x_\alpha : \alpha \in \Phi \}$ such that $G$ is generated by $\{ x_{\alpha}(t) : \alpha \in \Phi, t \in k \}$ modulo the relations
\begin{itemize}
\item[(A)] $x_{\alpha}(s+t) = x_{\alpha}(s) x_{\alpha}(t)$. 

\item[(B)] $[x_{\alpha}(t), x_{\beta}(s) ] = \prod_{i,j}  x_{i\alpha +j \beta}(n_{\alpha, \beta, i,j} t^i s^j)$ 
where the product is taken over $(i,j)$ in lexicographic order and the  $n_{\alpha, \beta,i,j}$'s are certain integers which we will not make precise 
here.

\item[(B')] $w_{\alpha}(t) x_{\alpha}(s) w_{\alpha}(t)^{-1} = x_{-\alpha}(-st^{-2})$. 

\item[(C)] $h_{\alpha}(st) = h_{\alpha}(s)h_{\alpha}(t)$. 
\end{itemize}

\noindent The condition (B') is not necessary if $G \ne \upr{SL}_2$.  The elements $h_{\alpha}(t)$ for $\alpha \in \Delta$ generate the group $T$, and the elements $x_\alpha(t)$ for $\alpha \in \Phi^+$, $t \in k$, generate the group $U^+$.  Similarly, $\{ x_\alpha(t) : \alpha \in \Phi^-, t \in k \}$ generates the opposite unipotent $U^- = \alg{U}^-(k)$.      

Steinberg considered the group $\tilde{G}$ generated by elements $\tilde{x}_{\alpha}(t)$, for $\alpha \in \Phi$, modulo the (analogous) relations (A), (B) and (B') above (ignoring (C)). There is clearly a natural surjection 
\[  p_G: \tilde{G} \longrightarrow G \]
given by  $\tilde{x}_{\alpha}(t) \mapsto x_{\alpha}(t)$.  

From relations (A) and (B), it can be seen that $x_\alpha(t) \mapsto \tilde x_\alpha(t)$ extends to a homomorphism $\sigma^\pm \From U^\pm \Into \tilde G$ splitting $p_G$, i.e., $p_G \circ \sigma^\pm = \Id$.

Define elements $\tilde{w}_{\alpha}(t)$ and $\tilde{h}_{\alpha}(t)$ analagously to $w_\alpha(t)$ and $h_\alpha(t)$ above. We let 
$\tilde{T}$ denote the subgroup of $\tilde{G}$ generated by the $\tilde{h}_{\alpha}(t)$'s. Steinberg showed (\cite[Chapter 7, Theorem 10]{S2}):


\begin{thm} \label{T:Steinberg}
The group $\Ker(p_G)$ is central in $\tilde{G}$, and 
$$1 \To \Ker(p_G) \To \tilde G \To G \To 1$$
is a universal central extension of $G$.  In particular $\pi_1(G) = \Ker(p_G)$. Furthermore, $\pi_1(G) \subset \tilde{T}$.
\end{thm}


\subsection{ Steinberg's cocycles.}
In \cite[p.194]{Mo}, Moore described a 2-cocycle which represents the universal extension $\tilde{G}$, depending on choices of Weyl representatives and an ordering of the simple roots.  Let $\tilde{N}$ (respectively $N$) denote the subgroup of $\tilde{G}$ (resp.~$G$) generated by the $\tilde{w}_{\alpha}(t)$'s 
(resp.~$w_{\alpha}(t)$). Then 
\[  \tilde{N}/\tilde{T} \cong N/T  = W \]
is the Weyl group of $G$.
For each $w \in W$, we fix a representative $\tilde{w} \in \tilde{N}$ and denote its projection $p_G(\tilde{w}) \in N$ by $\dot w$.  

Each element $g \in G$ lies in a unique Bruhat cell $BwB$, and can be uniquely represented as:
\[  g = u_w \cdot w \cdot t \cdot u, \quad t = \prod_{\alpha \in \Delta} h_{\alpha}(t_{\alpha}) \in T, \, u \in U^+, \, u_w \in U_w = \prod_{\substack{\alpha > 0 \\ w\alpha < 0} } U_{\alpha}. \]
We define a section $s: G \longrightarrow \tilde{G}$ by setting
\[ s(g)  = \tilde{u}_w \cdot \tilde{w} \cdot \tilde{h} \cdot \tilde{u} \]
where $\tilde{u}_w = \sigma^+(u_w)$ and $\tilde{u} = \sigma^+(u)$ and $\tilde{h} = \prod_{\alpha \in \Delta} \tilde{h}_{\alpha}(t_{\alpha})$.  Here we must fix an ordering of the simple roots, in order for this product to make sense.

This gives a 2-cocycle $b_{\univ} \From G \times G \longrightarrow \pi_1(G)$, given by
\[  b_{\univ}(g_1, g_2) = s(g_1)s(g_2) s(g_1g_2)^{-1}. \]
We call this the {\bf universal Steinberg cocycle of $G$}. 

Given any central extension $E \in \CExt(G,A)$, there is a unique homomorphism $f \From \tilde G \To E$ lying over the identity map on $G$.  Defining $s_E = f \circ s$ we have a section $s_E \From G \To E$, and from this a cocycle
\[  c_E (g_1, g_2) = f(b_{\univ}(g_1, g_2))  \in A \]
which incarnates $E$.  We call $c_E$ the {\bf Steinberg cocycle of $E$}. 


\subsection{ Moore's upper bound for $\pi_1(G)$.}
After Steinberg's construction of the universal central extension by generators and relations, 
one can hope for an explicit presentation of $\pi_1(G)$. This question was taken up in \cite{Mo}.
His results are summarized in the following theorem (cf. \cite[Lemma 8.1, Theorem 8.1,
Lemma 8.2 and Lemma 8.4]{Mo}).

\begin{thm} \label{T:Moore}
(i) $\pi_1(G)$ is generated by the elements $b_{\univ}(h_{\alpha}(s), h_{\alpha}(t))$ 
for $\alpha \in \Delta$ and $s, t \in k^{\times}$.
In fact, if we fix a long root $\alpha_0$, then $\pi_1(G)$ is generated by the elements 
$b_{\univ}(h_{\alpha_0}(s), h_{\alpha_0}(t))$ for $s, t \in k^{\times}$.  

(ii) If $c = f \circ b_{\univ}$ is a Steinberg cocycle valued in $A$ (with $f \From \pi_1(G) \rightarrow A$), then 
$c$ is completely determined by its restriction to $T \times T$. In fact, if $\alpha$ is a long root and $T_{\alpha}$ the 1-dimensional torus 
generated by $h_{\alpha}(t)$, then $c$ is completely determined by its restriction to $T_{\alpha} \times T_{\alpha}$.
\end{thm}


\begin{cor} \label{C:Moore}
Let $\alpha$ be a long root in $\Phi$ and let $G_{\alpha} \cong \upr{SL}_2$ be the subgroup generated by $U_{\alpha}$ and $U_{-\alpha}$. 
Then for any $A$, the natural map
\[   H^2(G, A) \longrightarrow H^2(G_{\alpha}, A) \]
is an injection. Equivalently, the natural map
\[  \pi_1(G_{\alpha}) \longrightarrow \pi_1(G) \]
is surjective.
\end{cor}

\noindent This corollary is the key tool in the analysis of extensions of a split group.

The theorem significantly reduces the number of generators needed for $\pi_1(G)$. Indeed, 
in view of the theorem, we fix a long root $\alpha$ and set
\[  b_{\univ, \alpha}(s,t) = b_{\univ}(h_{\alpha}(s), h_{\alpha}(t)) \]  
so that 
\[  b_{\univ, \alpha}: k^{\times} \times k^{\times} \longrightarrow \pi_1(G). \]
We know that $\pi_1(G)$ is generated by $b_{\univ, \alpha}(s,t)$ for $s, t \in k^{\times}$. Moreover, Moore showed that 
under a simple condition, the function $b_{\univ, \alpha}$ is bimultiplicative:

\begin{prop}
If there exists a root $\beta$ such that $\langle \beta^{\vee}, \alpha \rangle =1$, then $b_{\univ,\alpha}$ is bimultiplicative.
This condition holds as long as $G \ne \upr{Sp}_{2n}$ ($n \geq 1$). If $G = \upr{Sp}_{2n}$ (e.g., if $G = \upr{SL}_2$), then $b_{\univ,\alpha}$
may not be bimultiplicative.
\end{prop}

Now we want to know what are the relations satisfied by the $b_{\univ, \alpha}(s,t)$.
By working explicitly with the group $\upr{SL}_2$, Moore was able to show (cf. \cite[Theorem 9.2]{Mo}):

\begin{thm} \label{T:Moore2}
If $G = \upr{SL}_2$, then $\pi_1(G)$ is the group generated by the $b(s,t) := b_{\univ, \alpha}(s,t)$ subject to the relations:
\begin{itemize}
\item[(1)] (normalized cocycle identities) 
\[  b(st,r)b(s,t) = b(s,tr)b(t,r), \quad b(s,1) = b(1,s) = 1. \]

\item[(2)] $b(s, t) = b(t^{-1},s)$.

\item[(3)] $b(s,t) = b(s,-st)$. 

\item[(4)] $b(s,t) = b(s, (1-s)t)$ if $s \ne 1$.
\end{itemize}
There is in fact some redundancy in these relations: under (1) and (4), (2) and (3) are equivalent.
\end{thm}

\subsection{ Definition:} We call functions $c: k^{\times} \times k^{\times} \longrightarrow A$ satisfying the above identities (1)-(4) 
{\bf $A$-valued Steinberg cocycles on $k^{\times}$}. We denote the set of such functions by $\St(k^\times, A)$.  Since the universal cocycle $b_{\univ,\alpha}$ is bimultiplicative if $G \ne \upr{Sp}_{2n}(k)$, we define a subgroup $\St^\circ(k^\times, A) \subset \St(k^\times, A)$ consisting of those $A$-valued Steinberg cocycles which are bimultiplicative. 

The elements of $\St^\circ(k^\times, A)$ can be more simply described as those maps $c \From k^\times \times k^\times \To A$ satisfying just two conditions:

\begin{itemize}
\item[(1')] (bimultiplicative) $c(rs, t) = c(r,t) \cdot c(s,t)$ and $c(r,st) = c(r,s) \cdot c(s,t)$. 

\item[(2')]  $c(s, 1-s) = 1$ if $s \ne 1$. 
\end{itemize}

\noindent The relations (1') and (2') are important in algebraic K-theory.  Namely, they occur in the definition of the {\bf Milnor-Quillen $\alg{K}_2$-group} \cite{Mi}.  This is the abelian group  
\[  \alg{K}_2(k) = \frac{k^{\times} \otimes_{\ZZ} k^{\times} }{\langle x \otimes (1-x) : x \neq 1 \rangle}. \]


Thus $\St^\circ(k^\times, A) = \Hom( \alg{K}_2(k), A)$.   A corollary of the above discussions and Pontryagin duality is:

\begin{cor} \label{C:Moore2}
(i) For $G \ne \upr{Sp}_{2n}(k)$ (resp.~$G = \upr{Sp}_{2n}(k)$) and any $A$, there is a natural inclusion
\[  H^2(G,A) \hookrightarrow  \St^\circ(k^\times,A) \quad \text{(resp. $\St(k^\times, A)$).} \]


(ii) If $G = \upr{Sp}_{2n}(k)$, then $\pi_1(G)$ is a quotient of the 
group generated by $b_{\univ, \alpha}(s,t)$ subject to the relations (1) - (4) of Theorem 
\ref{T:Moore2}. 

(iii) If $G \ne \upr{Sp}_{2n}(k)$, then $\pi_1(G)$ is a quotient of the  
group $\alg{K}_2(k)$. 
\end{cor}

However, for a general Chevalley group $G$, Moore was unable to determine  
whether the relations in the corollary are enough to define $\pi_1(G)$, or whether more relations are necessary. 

\subsection{ Matsumoto's determination of $\pi_1(G)$.}
In \cite{Ma}, Matsumoto was able to complete Moore's results by showing:

\begin{thm} \label{T:Mat}
(ii) If $G = \upr{Sp}_{2n}(k)$, then $\pi_1(G)$ is isomorphic to the 
group generated by $b_{\univ, \alpha}(s,t)$ subject to the relations (1) - (4) of Theorem 
\ref{T:Moore2}. Thus
\[  H^2(G,A) \cong \St(k^\times, A). \]


(iii) If $G \ne \upr{Sp}_{2n}(k)$, then $\pi_1(G)$ is isomorphic to $\K_2(k)$.  Thus,
\[  H^2(G,A) \cong \St^\circ(k^\times,A). \]
\end{thm}

We remark that since we have an upper bound $H^2(G,A) \hookrightarrow \St(k^\times,A)$ (or $\St^\circ(k^\times, A)$) from Moore, to show that this upper bound is attained is a question of construction of central extensions.  Namely, given an element of $f \in \St(k^\times,A)$ or $\St^\circ(k^\times, A)$, one needs to construct a central extension of $G$ by $A$ whose associated Steinberg cocycle gives rise to $f$. This was what Matsumoto did.

\section{ Groups over Local Fields}

In this section, we consider the main problem highlighted in the first section for groups over local fields. 
Let $k$ be a local field and let $\alg{G}$ be a (algebraically) simply-connected semisimple group over $k$. 
We set $G = \alg{G}(k)$, so that $G$ is a topological group and we are interested in its topological central extensions.

If $\alg{G}$ is $k$-isotropic, it is known that $G$ is topologically perfect, so that there is a chance that $\pi_1(G)$ 
exists.  The main result we want to highlight here is:

\begin{thm} \label{T:local}
Suppose that $k$ is non-archimedean.
Assume that $\alg{G}$ is absolutely simple and $k$-isotropic. Then
\[ H^2(G, S^1) \cong \mu(k)^{\vee}, \]
where $\mu(k)$ denotes the finite group of roots of unity contained in $k$. 
In particular, $\pi_1(G)$ exists and is equal to $\mu(k)$.
\end{thm}

We make several remarks:

\noindent (1) The assumption that $k$ is non-archimedean is for convenience: it allows us to give a simple statement. Over $\mathbb{R}$ or 
$\mathbb{C}$, the situation is completely understood.

\noindent (2) The condition that $\alg{G}$ be absolutely simple is not crucial. If $\alg{G}$ is just semisimple and simply-connected, then
\[  \alg{G} \cong \prod_i  \text{Res}_{k_i/k} \alg{G}_i \]
with $\alg{G}_i$ absolutely simple. Thus if each $\alg{G}_i$ is $k_i$-isotropic, the theorem implies that
\[  \pi_1(G) = \prod_i \mu(k_i).  \]

\noindent (3) The theorem is the culmination of the work of Steinberg \cite{S1}, Moore \cite{Mo}, Matsumoto \cite{Ma}, Deodhar \cite{De}, Prasad-Raghunathan \cite{PR1,PR2, PR3}, Prasad-Rapinchuk \cite{PR}, G. Prasad  \cite{P} and 
Deligne \cite{D2}. In the rest of the section, we will describe some ideas in its proof. 

\noindent (4) If $\alg{G}$ is anisotropic, absolutely simple, and simply-connected, then $G = \upr{SL}_1(D)$ where $D$ is a division algebra over $k$.  In this case, one can still demand to compute $H^2(G, S^1)$, even though $G$ is not perfect. Such a computation was done by Prasad-Rapinchuk. We will not discuss this here.


\subsection{ The case of split groups.}
When $\alg{G}$ is split, the theorem was proved by the combined work 
of Moore and Matsumoto, which made decisive use of the analysis of the 
abstract universal central extension given in the last section. Let $G_{\abst}$ denote $\alg{G}(k)$ regarded as an abstract group (with discrete 
topology).  Then since any topological central extension is an abstract extension, we have a natural map
\[  H^2(G,A) \longrightarrow  H^2(G_{\abst}, A_{\abst}). \]
It turns out that this natural map is always an inclusion (for any topological group $G$), so that there is a natural surjection
\[  \pi_1(G_{\abst}) \longrightarrow \pi_1(G). \]
By Theorem \ref{T:Mat}, we know that
\[  H^2(G_{\abst}, A_{\abst}) \cong \St(k^\times, A) \text{ or } \St^\circ(k^\times, A). \]
Thus it remains to determine which $A$-valued Steinberg cocycles correspond to topological extensions. The following result is
both simple (to absorb) and natural:

\begin{thm}
Let $E \in \CExt(G_{\abst}, A_{\abst})$. Then the following are equivalent:

\begin{itemize}
\item[(i)]  The Steinberg cocycle $c_E: G \times G \longrightarrow A$ of $E$ is Borel measurable.

\item[(ii)] $c_E$ is continuous on $T \times T$.

\item[(iii)] $c_E$ is continuous on $T_{\alpha} \times T_{\alpha}$.  ($\alpha$ a long root as before.)

\item[(iv)] $E$ is a topological central extension.
\end{itemize}
\end{thm}  		

Thus, to classify topological central extension, we are reduced to classifying the set $\St_{\cont}(k^\times, A)$ of continuous $A$-valued Steinberg cocycles. This problem was solved by Moore (cf. \cite[Chapter 2]{Mo}). To describe his answer, we first recall that there is a natural supply of 
elements of $\St_{\cont}(k^\times, A)$ arising from local class field theory. Namely, 
if we let $\mu = \# \mu(k)$, there is a surjective $\mu$-power residue symbol
\[  (-,-): k^{\times} \times k^{\times} \longrightarrow \mu(k). \]
Moore observes \cite[Chapter II.(3)]{Mo} that $(-,-)$ is an element of $\St_{\cont}^\circ(k^\times, \mu(k))$.

Now given any $A$ and a homomorphism $f : \mu(k) \longrightarrow A$, we obtain an element
\[  f \circ (-,-) \in \St_{\cont}^\circ(k^\times, A) \subset \St_{\cont}(k^\times, A). \]
Thus we have a map
\[  \Hom(\mu(k), A) \longrightarrow \St_{\cont}(k^\times, A). \]
The result of Moore \cite[Theorem 3.1]{Mo} is:

\begin{thm}
The natural map above is bijective:
\[  
\St_{\cont}^\circ(k^\times, A) = \St_{\cont}(k^\times, A) = \Hom(\mu(k), A).
\]
In particular, each element of $\St_{\cont}(k^\times, A)$ is bimultiplicative (recall that $k$ is non-archimedean here) and $\pi_1(G) = \mu(k)$.
\end{thm}

\subsection{ Deodhar's work for quasi-split groups.}
Using a similar generators-relations approach based on a Chevalley-Steinberg system of \'{e}pinglage, Deodhar  \cite{De}
was able to extend Moore's results to the case when $\alg{G}$ is quasi-split. In particular, he obtained an upper bound for $\pi_1(G)$,
namely that
\[  \mu(k) \twoheadrightarrow \pi_1(G). \]
Once again, to establish that this is a bijection, one needs to construct topological central extensions. Thankfully, in this case, 
one does not need to give new constructions of central extensions.
One can finesse the difficulty by using an observation of Deligne (unpublished) to reduce to the case of split groups.

Deligne's observation makes use of one consequence of Matsumoto's work which is useful to know. Suppose we have an embedding
\[  i: \alg{SL}_2 \hookrightarrow \alg{G}. \]
Let $\alg{H}$ be a maximal split torus of $\alg{SL}_2$ and $\alg{T}$ a maximal split torus of $\alg{G}$ containing $i(\alg{H})$. 
There is then an embedding of $\ZZ$-modules:
\[  X_*(\alg{H}) \hookrightarrow X_*(\alg{T}). \]
Fix a Weyl group invariant inner product $\langle - ,- \rangle$ on $X_*(\alg{T}) \otimes \RR$ such that  for any long root $\alpha$ (so that 
$\alpha^{\vee}$ is short),
\[  \langle \alpha^{\vee}, \alpha^{\vee} \rangle  = 1. \]
Now take any generator $\mu$ of $X_*(\alg{H}) \cong \ZZ$ and set
\[   n(i, G) = \langle \mu, \mu \rangle \geq 1. \]
Consider the induced map
\[  i^*: H^2(G, A) \longrightarrow H^2(\upr{SL}_2, A). \]
Then the following lemma follows from [Mat, Lemma 5.4] and its proof.

\begin{lemma}
The image of $i^*$ is $n(i,G) \cdot H^2(\upr{SL}_2,A)$.
\end{lemma} 

Now Deligne showed that for each quasi-split $\alg{G}$, one can find
\begin{itemize}
\item a split group $\alg{G}'$ containing $\alg{G}$,
\item an embedding $i: \alg{SL}_2 \hookrightarrow \alg{G}$
\end{itemize}
such that $n(i, G') = 1$. Thus the composite
\[  H^2(G', A) \longrightarrow H^2(G,A)  \longrightarrow H^2(\upr{SL}_2, A) \]
is surjective. In particular, one can deduce that $H^2(G,A) \cong H^2(\upr{SL}_2, A) \cong \Hom(\mu(k), A)$.


\subsection{ The work of Prasad-Raghunathan for general $k$-isotropic groups.}

When $\alg{G}$ is $k$-isotropic but not quasi-split, then the above strategy is not feasible because we do not have an explicit description of 
$\pi_1(G_{\abst})$ to begin with. In this case, Prasad and Raghunathan \cite{PR1, PR2} have to resort to more geometric ideas (using the Bruhat-Tits building of $G$)
in order to compute $H^2(G, S^1)$.  The details are too intricate to discuss here. In the end, they 
showed that $\pi_1(G)$ is a quotient of $\mu(k)$ with kernel at most of size $2$.  This was then strengthened 
to an isomorphism using the results of \cite{D2} and \cite{PR}.


\section{ Adelic Groups}

In this section, suppose that $k$ is a global field and let $\mathbb{A}$ be its adele ring. For each place $v$ of $k$, let $k_v$ be the corresponding 
completion of $k$.  Let $\alg{G}$ be a simply-connected semisimple group over $k$. We set
\[  G_k = \alg{G}(k), \quad G_{\AA} = \alg{G}(\AA), \quad G_v = \alg{G}(k_v). \]
If $S$ is a finite set of places of $F$, one may also work with the $S$-adeles $\AA_S$; then $G_{\AA_S}$ is the restricted product of the $G_v$ for 
$v \notin S$.  There is a natural diagonal map $i: G_k \hookrightarrow G_{\AA_S}$, and one is interested in classifying topological central extensions of $G_{\AA}$ which split over $i(G_k)$.  These are classified by
\[  M(S, G) = \Ker(H^2(G_{\AA_S}, S^1) \longrightarrow H^2(G_k, S^1)). \]
This group is called the {\bf $S$-metaplectic kernel}. If $S = \emptyset$, we call it the {\bf absolute metaplectic kernel} and denote it simply by 
$M(G)$. The computation of $M(S,G)$ was achieved after a long series of papers by Prasad-Raghunathan \cite{PR1, PR2, PR3} and Prasad-Rapinchuk \cite{PR}. 


One reason for focusing on central extensions of $G_{\AA}$ which become split over $G_k$ is that one is eventually interested in the theory of automorphic forms of coverings $\tilde{G}_{\AA}$ of $G_{\AA}$: these are functions on $i(G_k) \backslash \tilde{G}_{\AA}$.  Another reason is that the computation of $M(S,G)$ arises in the study of the congruence subgroup problem. 

\subsection{ Local-to-global.} 
If $\alg{G}$ is $k$-isotropic, then $\pi_1(G_v)$ exists for all $v$ by Theorem \ref{T:local}. Moreover, it is known that if $K_v$ is a hyperspecial maximal compact subgroup of $G_v$, then $K_v$ is perfect  , so that $\pi_1(K_v)$ exists for almost all $v$ by Proposition \ref{P:suff}(iii). Moreover, for almost all $v$, 
the natural map
\[  i_{v*}: \pi_1(K_v) \longrightarrow \pi_1(G_v) \]
is the zero map. Thus by Proposition \ref{P:product}, $\pi_1(G_{\AA_S})$ exists and is equal to
\[  \pi_1(G_{\AA_S}) = \bigoplus_{v \notin S} \pi_1(G_v). \]
Further, the discrete group $G_k$ is perfect so that $\pi_1(G_k)$ exists also. Thus, using Proposition \ref{P:relative},
we have  the relative fundamental group
\[  \pi_1(G_{\AA_S},  G_k) = \left( \bigoplus_{v \notin S} \pi_1(G_v) \right) /  \overline{i_*(\pi_1(G_k))}. \]
Given all these, one deduces that the functor 
\[ A \mapsto M(S, G, A) = \Ker(H^2(G_{\AA_S}, A) \longrightarrow H^2(G_k, A)) \]
is represented by $\pi_1(G_{\AA_S}, G_k)$, so that 
\[  M(S,G) = \Hom(\pi_1(G_{\AA_S}, G_k), S^1). \]
 

Thus, if $\alg{G}$ is $k$-isotropic, the problem of computing $M(S,G)$ is the same as computing the relative fundamental group 
$\pi_1(G_{\AA_S}, G_k)$. Since we know the local $\pi_1(G_v)$'s very explicitly, one approach to computing $\pi_1(G_{\AA_S}, G_k)$ is 
to describe as explicitly as possible the closure of the image of $\pi_1(G_k)$. For this, one would need to know $\pi_1(G_k)$ very explicitly. As we 
noted in Section 2, we have this explicit description when $\alg{G}$ is quasi-split, thanks to the work of Steinberg, Moore, Matsumoto 
and Deodhar. When $\alg{G}$ is not quasi-split, such an approach to computing $\pi_1(G_{\AA_S}, G_k)$ is not feasible. This is why for 
non-quasi-split groups, Prasad-Raghunathan and Prasad-Rapinchuk have to resort to completely different ideas to solve this problem.

In any case, the main global theorem of \cite{PR} is:

\begin{thm}
Let $\alg{G}$ be an absolutely simple, simply connected semisimple group over $k$.  If $\alg{G}$ is a special unitary group over a noncommutative division algebra, assume a certain conjecture (U).  Let $S$ be a finite (possibly empty) set of places of $k$.  Then we have:
\begin{enumerate}
\item[(i)] $M(S,G) \subset \mu(k)^{\vee}$.

\item[(ii)] If $S$ contains a non-archimedean place $v$ where $G_v$ is isotropic or a real place where $G_v$ is not topologically simply-connected,
then $M(S,G) = 0$.

\item[(iii)] If $S = \emptyset$, then $M(S,G) = \mu(k)^{\vee}$.
\end{enumerate}
\end{thm}


As for the local theorem, this theorem is the culmination of the work of many people, culminating in the eventual work of Prasad-Rapinchuk \cite{PR}.

 \section{ Brylinski-Deligne Theory}  \label{S:BD}
 As our brief discussion of the historical development of the structure theory of covering groups shows, much of the earlier work is focused on determining the fundamental group or the universal central extension. This almost immediately restricts one to the case when $\alg{G}$ is a simply-connected linear algebraic group over a field $k$. 
 The disadvantage of this is that it is a common strategy in Lie theory to prove results by induction through Levi subgroups of parabolic subgroups. However, the Levi subgroups are only reductive groups and not semisimple. Thus the structure theory of (topological) central extensions obtained in previous sections does not apply to the Levi subgroups of $\alg{G}$.
 
 
 When $\alg{G}$ is not simply-connected, for example if $\alg{G}$ is a special orthogonal group, a typical way of obtaining central extensions of $G = \alg{G}(k)$ is to fix an embedding $\upr{G} \hookrightarrow \upr{SL}_r$ and then to pullback some known central extensions of $\upr{SL}_r$. For example, one may embed $\upr{GL}_r$ into $\upr{SL}_{r+1}$ as an $r \times r$ block, and pullback a (topological) central extension of $\upr{SL}_{r+1}$; this gives an extension which has been studied in some detail by Kazhdan and Patterson \cite{KP1} and is one member of a family of covers of $\upr{GL}_r$ known as the Kazhdan-Patterson covers. 
 
While such constructions give examples of covering groups, with some control on their structure through one's knowledge of the relevant 2-cocycles on $\upr{SL}_r$, they do not amount to a systematic theory or classification.
 

In their 2001 IHES paper \cite{BD}, Brylinski and Deligne  approached the subject from a different angle. They returned to the very neat results obtained in the split simply-connected case by Steinberg, Moore and Matsumoto, where one has an extension of abstract groups 
\[  \begin{CD}
1 @>>> \alg{K}_2(k) @>>>  \tilde{G}  @>>>  G= \alg{G}(k) @>>>1,  \end{CD} \]
which is universal if $\alg{G}$ is not of type $C$. Their idea (from our perspective) is to ``remove the $k$" in the above short exact sequence. More precisely, regarding $\alg{K}_2$ and $\alg{G}$ as sheaves of groups on the big Zariski site of ${\rm Spec}(k)$,  they consider the problem of understanding or classifying the central extensions of group sheaves
\[  \begin{CD}
 1@>>>  \alg{K}_2 @>>> \tilde{\alg{G}}  @>>> \alg{G} @>>> 1. \end{CD} \]
Such a $\tilde{\alg{G}}$ is also called a multiplicative $\alg{K}_2$-torsor over $\alg{G}$ and the problem is to give a classification of the Picard category of such multiplicative $\alg{K}_2$-torsors with $\alg{G}$ fixed, i.e., to describe this category in simpler terms. Brylinski-Deligne managed to give a very reasonable answer to this classification problem which depends functorially on $\alg{G}$. Their results  will be summarised and described in the papers in this volume. 

Suppose one has a multiplicative $\alg{K}_2$-torsor $\tilde{\alg{G}}$ over a local field $k$. Then on taking $k$-points, one obtains a central extension of discrete groups
\[  \begin{CD}  
1 @>>>  \alg{K}_2(k) @>>>  \tilde{\alg{G}}(k) @>>> \alg{G}(k)  = G @>>> 1. \end{CD} \]
Here the sequence remains exact on the right because $H_{\Zar}^1(k, \alg{K}_2)  = 1$.  If one pushes this sequence out via the norm residue symbol $\K_2(k) \longrightarrow \mu(k)$, then one obtains a topological central extension
\[ \begin{CD} 1 @>>>  \mu(k) @>>> \tilde{G} @>>> G @>>> 1. \end{CD} \]
Thus, multiplicative $\alg{K}_2$-torsors over local fields give rise to topological central extensions. Such topological central extensions are thus of ``algebraic origin". 

Now suppose $k$ is a global field with ring of adeles $\AA$. Then the analog of the above construction shows that one inherits a central extension
\[  \begin{CD}
 1@>>> \mu(k) @>>>  \tilde{G}_{\AA}  @>>>  G_{\AA} @>>> 1  \end{CD} \]
 from a multiplicative $\alg{K}_2$-torsor over $k$. A key feature of the Brylinski-Deligne theory is that this central extension of $G_{\AA}$ comes equipped with a canonical splitting $G_k \hookrightarrow \tilde{G}_{\AA}$; this follows from reciprocity for norm residue symbols. In other words, a multiplicative $\alg{K}_2$-torsor over a global field gives rise to a topological central extension of the adelic group $G_{\AA}$
together with a splitting over $G_k$.  This means that one is immediately in a position to begin the study of automorphic forms of $\tilde{G}_{\AA}$. 


\section{ Representation Theory and Automorphic Forms}
In this section, we shall give a brief discussion of several representative works on the representation theory, harmonic analysis and the theory of automorphic forms on covering groups.  

Shortly after the middle of the last century, the classical theory of integer weight modular forms on the upper half plane was recast in the framework of automorphic representations on the group $\alg{SL}_2$. On the other hand, it has been known that fractional weights modular forms exist and play a significant role in classical modular form theory. One early example is the Jacobi theta function (a weight 1/2 modular form), or more generally the theta function associated to an integer lattice of odd rank.  It was then observed that such modular forms should correspond to automorphic representations on covering groups of $\alg{SL}_2(\AA)$. This gives a strong impetus for a systematic study of covering groups of adelic groups.


\subsection{ Segal-Shale-Weil representation.}  One of the first systematic study of representations of a covering group is the work of Weil \cite{We1} on the so-called Segal-Shale-Weil representations (also called the oscillator representations) of the unique 2-fold cover of $\upr{Sp}_{2n}(k)$ (where $k$ is a local field). This 2-fold cover is called the metaplectic group $\upr{Mp}_{2n}(k)$.  Weil and others after him (such as Kubota \cite{Ku1} and Rao \cite{Ra}) gave a comprehensive study of the 2-cocycles describing the metaplectic groups and its Weil representations. Weil's goal for developing this was to  reformulate the theory of theta functions in the representation theoretic framework and to express previous results of Siegel (such as the Siegel mass formula and  Siegel-Weil formula) in this language \cite{We2}. The Weil representations subsequently became a key ingredient in Howe's theory of dual pair correspondence (or theta correspondence).  

\subsection{ The work of Kubota and Patterson.}  \label{SS:kubota}
In the late 1960's, almost concurrently as Moore and Matsumoto were doing their groundbreaking work, Kubota initiated a systematic study of  the coverings of $\upr{SL}_2$ or $\upr{GL}_2$ (beyond the 2-fold cover), giving precise 2-cocycles for these covers \cite{Ku1, Ku2}. He was also interested in constructing analogs of Jacobi's theta function for these higher degree covers, using the residues  of Eisenstein series.  Patterson made a detailed study of the Fourier expansion of some of 
these theta functions on higher degree covers, noting that they contain interesting arithmetic information. In particular, for the 3-fold cover of $\upr{GL}_2$, he showed in \cite{P1,P2} that the Fourier coefficients of the cubic theta function are  cubic Gauss sums. Using this connection, Heath-Brown and Patterson \cite{HP} showed the equidistribution of the angular components of cubic Gauss sums.  This suggests that one might find arithmetic applications by studying the Fourier expansion of interesting automorphic forms on covering groups. We will discuss some other of these arithmetic applications later on. 

For higher degree covers, however, the structure of the Fourier coefficients of the generalized theta functions becomes much more complicated. This was subsequently explained by Deligne \cite{D1}  as a consequence of the fact that Whittaker models are not unique for higher degree covers of $\alg{SL}_2$.


\subsection{ Shimura's correspondence.}   \label{SS:shimura}
One of the key milestones in the theory of automorphic forms on covering groups is Shimura's 1973 Annals paper \cite{Sh}, in which he developed a theory of Hecke operators for half integer weight modular forms and proved a correspondence between half integer weight modular forms and integer weight modular forms. Shimura proved the correspondence which bears his name by using the converse theorem of Weil. To do so, he introduced  another innovation in his paper: a Rankin-Selberg  integral for the standard L-function of a half-integer weight modular form. A slight variant of this Rankin-Selberg integral gives the symmetric square L-function of an integer weight modular form, which was used by Gelbart-Jacquet in their work on the symmetric square lifting from $\upr{GL}_2$ to $\upr{GL}_3$. 


This influential paper of Shimura is the first to establish a lifting from Hecke eigenforms of a covering group to those of a linear group. It led to two independent lines of development, as we recall below. Both of these arise from the attempt to formulate Shimura's results in the setting of automorphic representations.


\subsection{ Kazhdan-Patterson covering and  Flicker-Kazhdan lifting.}  \label{SS:KP}
The first line of development from Shimura's paper is the work \cite{F} of Flicker, who used the trace formula approach to prove the Shimura correspondence. More precisely, Flicker compared the trace formula of a particular degree $n$ cover of $\upr{GL}_2$ constructed by Kubota with that of $\upr{GL}_2$ and proved a one-to-one correspondence between cuspidal automorphic representaiotns of $\widetilde{\upr{GL}}_2$ and cuspidal representations of $\upr{GL}_2$ whose central character is an $n^{\th}$ power.  Thus,  his work went beyond what Shimura did as he considered not just 2-fold covers of $\upr{GL}_2$. In this adelic treatment of the Shimura correspondence, there is a local correspondence between genuine representations of the local covering group and those of the linear group $\upr{GL}_2$. This local correspondence is expressed by a local character identity. 

 Following up on this work, Kazhdan-Patterson \cite{KP1} considered degree $n$ covers of $\upr{GL}_r$ which are obtained by pulling back from
the degree $n$ cover of $\upr{SL}_{r+1}$ (with $\upr{GL}_r$ embedded in $\upr{SL}_{r+1}$ in a standard way) and a standard twisting operation. 
Such covers are now called Kazhdan-Patterson covers and they generalise the Kubota covers of $\upr{GL}_2$. 
In \cite{KP1}, Kazhdan and Patterson were largely interested in extending the results of Kubota noted above to the higher rank case of $\upr{GL}_r$; in particular, they constructed generalisations of theta functions as residues of Eisenstein series.   Further, their paper also laid the groundwork for an extension of Flicker's results from covers of $\upr{GL}_2$ to the Kazhdan-Patterson covers of $\upr{GL}_r$. This extension was initiated in their subsequent paper \cite{KP2} and pursued further in the paper \cite{FK} of Flicker-Kazhdan. In \cite{FK}, the authors  used a simple trace formula to prove a correspondence between cuspidal automorphic representations of the Kazhdan-Patterson covers of $\upr{GL}_r$ and cuspidal representations of $\upr{GL}_r$ under some simplifying local conditions (which allow the use of the simple trace formula).  As in the $\upr{GL}_2$ case, subordinate to this global correspondence is a local correspondence of representations based on a local character identity.

Somewhat unfortunately, it has been noted by several people that there are some errors in the papers \cite{KP1} and  \cite{FK}; given the subtlety of the structure theory of covering groups, this is quite understandable and certainly does not detract from the pioneering nature of these papers.  
The authors of \cite{FK} and \cite{KP1} have however not  provided an account and erratum for these errors. This is quite unfortunate, as there is no doubt that most of the results there must be true, at least if one imposes some conditions on the degree of the cover. Some further work in this direction, which cleared up some of the issues,  were carried out by A. Kable \cite{Ka}, P. Mezo \cite{Me1, Me2} and by Banks-Levy-Sepanski \cite{BLS} among others. 

\subsection{ The work of Waldspurger.}  \label{SS:walk}
Another line of development originating from Shimura's 1973 paper is the work of Waldspurger which uses the technique of theta correspondence. In two papers \cite{Wa1, Wa3}, using the theta correspondence for $\upr{Mp}_2 \times \upr{SO}_3$, Waldspurger obtained a complete description of the automorphic discrete spectrum of the metaplectic group $\upr{Mp}_2$ in terms of that of $\upr{PGL}_2 = \upr{SO}_3$. In particular, over local fields, Waldpsurger  obtained a classification of the genuine representations of $\upr{Mp}_2$ in the style of the local Langlands correspondence. Moreover, his description of the automorphic discrete spectrum of $\upr{Mp}_2$ is in the style of the Arthur conjecture \cite{Ar}, i.e.~using local and global packets and having a global multiplicity formula. 

What is especially intriguing is that the global multiplicity formula involves the global root number of a cuspidal representation of $\upr{PGL}_2$.  It should be noted, however, that at a critical point of his work, Waldspurger had needed to appeal to Flicker's results \cite{F} obtained by the trace formula mentioned above. As a consequence, Waldspurger showed the existence of nonvanishing of central L-value of quadratic twists of automorphic L-functions of $\upr{GL}_2$. 
Nowadays, however, one can avoid appealing to \cite{F} , as Bump, Friedberg and Hoffstein have independently proven the necessary nonvanishing of central L-values. For a more detailed discussion of Waldspurger's work in light of Bump-Friedberg-Hoffstein \cite{BFH1, BFH2, FH}, the reader can look at \cite{G1}.

The work of Waldspurger has led to other significant arithmetic applications. For example, in \cite{Wa2}, he obtained a formula expressing the Fourier coefficients of half integral weight modular forms in terms of the central critical L-value of its Shimura correspondent. This result was applied by Tunnell \cite{Tu} to provide a solution to the congruence number problem modulo the Birch-Swinnerton-Dyer conjecture. 


It is natural to ask for an extension of Waldspurger's work to $\upr{Mp}_{2n}$ by using the theta correspondence for $\upr{Mp}_{2n} \times \upr{SO}_{2n+1}$. This has come slowly over the past 35 years. In the local setting, the local Shimura correspondence, giving a classification of the genuine representations of $\upr{Mp}_{2n}$ in terms of the representations of $\upr{SO}_{2n+1}$, was shown by Adams-Barbasch  \cite{AB} over $\RR$ in the 1990's.  The analogous result for $p$-adic fields was only shown fairly recently by Gan-Savin \cite{GS}. 
Over number fields, a conjectural extension of Waldspurger's results to $\upr{Mp}_{2n}$ was given in \cite{GGP} and \cite{G2}.
A recent preprint of Gan-Ichino gives a classification of the part of the automorphic discrete spectrum of $\upr{Mp}_{2n}$ associated to tempered A-parameters, thus proving the conjecture formulated in \cite{GGP}.  The reason for this long lapse in extending Waldspurger's results to $\upr{Mp}_{2n}$ is because that it requires recent advances in the theory of theta correspondences as well as the recent results of Arthur \cite{Ar} on the automorphic discrete spectrum of classical groups. The result of Waldspurger on Fourier coefficients of half-integral weight modular form was extended to the setting of  the Whittaker-Fourier coefficients of cuspidal representations of $\upr{Mp}_{2n}$ in a recent series of papers by Lapid-Mao \cite{LM1, LM2, LM3}.


 \subsection{ Fourier coefficients of metaplectic Eisenstein series and generalized theta functions.}
 The work of Kubota and Patterson on the Fourier coefficients of generalized theta functions and metaplectic Eisenstein series was continued by several mathematicians in the 1980's and 1990's, most notably Bump, Friedberg, Hoffstein and their students or collaborators. An early work   is the paper \cite{BH1} of Bump-Hoffstein which shows that cubic L-functions  occur in the Fourier expansion of Eisenstein series on a 3-fold Kazhdan-Patterson cover of $\upr{GL}_3$. Several conjectures were highlighted and formulated in the paper \cite{BH4} of Bump-Hoffstein and some further works in this direction include \cite{BH3, BBL, Su1, Su2}. There are complementary local results \cite{BFH3} on Whittaker functions of unramified genuine representations, analogous to the Casselman-Shalika formula in the linear case.
  These results on Fourier coefficients have  found many stunning arithmetic applications, concerning nonvanishing of central values of twists of automorphic L-functions, 
  such as \cite{BFH1, BFH2, BrFH}. 

  
  These early works ultimately led Brubaker,  Bump, Friedberg and Hoffstein to develop the theory of Weyl group multiple Dirichlet series in a series of papers (see for example, \cite{BBF, BBFH}), together with important contributions from Chinta and Gunnells \cite{CG1, CG2}.   This theory of Weyl group multiple Dirichlet series has found surprising connections with combinatorics, statistical physics and quantum groups.  
  The local theory of the Casselman-Shalika formula culminated in the recent papers of Chinta-Offen and  McNamara  \cite{CO, Mc1, Mc3}. 

Automorphic descent in the covering setting is also explored by Friedberg and Ginzburg \cite{FG2, FG3}.  Some other recent work on the Fourier coefficients of metaplectic Eisenstein series is contained in \cite{BrF, FG1, FZ}, and also the thesis work of Y.Q. Cai \cite{C1, C2}.


\subsection{ Automorphic L-functions}
Another active area of research concerns the theory of automorphic L-functions for metaplectic forms.
To have a definition of automorphic L-functions, one would need to have an understanding of unramified representations and a notion of the dual group of a covering group. 
For $\upr{Mp}_{2n}$ and the Kazhdan-Patterson covers of $\upr{GL}_n$ (under some conditions on the degree of covering), one has a natural candidate for the dual group. For $\upr{Mp}_{2n}$ the natural dual group is $\upr{Sp}_{2n}(\CC)$ and for the Kazhdan-Patterson covers, it is $\upr{GL}_n(\CC)$,  at least under some assumptions on the degree of cover.  For these examples of covering groups, one can define Satake parameters for unramified representations and thus define the notion of partial automorphic L-functions. To show the usual analytic properties of these L-functions,  one would try to find some Rankin-Selberg integrals for these automorphic L-functions. Some early work in this direction include \cite{BH2, BH3, BF}. 
On the other hand, the thesis work of D. Szpruch \cite{Sz1, Sz2} develops the Langlands-Shahidi theory for $\upr{Mp}_{2n}$.
A more recent preprint of Cai-Friedberg-Ginzburg-Kaplan \cite{CFGK} gives a sketch of a generalisation of the  doubling method of Piatetski-Shapiro-Rallis, which gives a Rankin-Selberg integral for the standard L-function of covers of classical groups.


Metaplectic forms have also proved useful in constructing Rankin-Selberg integrals for automorphic L-functions of linear groups. The prime
 example is the work of Bump-Ginzburg \cite{BG} which extended Shimura's original work to give a Rankin-Selberg integral for the symmetric square L-function of cuspidal representations of $\upr{GL}_n$, using an Eisenstein series on a double cover of $\upr{GL}_n$. Based on their work, the case of twisted symmetric square L-function is treated by Takeda \cite{Ta}.

\subsection{ Savin's Hecke algebra correspondence.}  \label{SS:savin}
In another direction,  Savin studied and determined the structure of the Iwahori Hecke algebra for covers $\tilde{G}$ of simply-connected groups  $G$ \cite{Sa1,Sa2} and showed that they are  isomorphic to the Iwahori Hecke algebra of an appropriate linear group. This gives a bijection between the irreducible genuine representations of $\tilde{G}$ with Iwahori-fixed vectors and those of the linear group. He did the same for the spherical Hecke algebra, thus obtaining a correspondence of unramified representations. These papers of Savin were the first to attempt a systematic development of the representation theory of general covering groups, going beyond treating families of examples.  It
 gives strong suggestions for the dual groups of covers of simply-connected groups. 


\subsection{ Character identitites.} 
The Flicker-Kazhdan local correspondence suggested that lifting of representations between covering and linear groups can be formulated in terms of local character identities. Such local character identities, in the context of Kazhdan-Patterson covers and other covering groups such as $\upr{Mp}_{2n}$, were studied by J.~Adams in a series of papers in the 1990's \cite{A1, A2}. Most of Adams' work is focused on coverings of real groups. It culminates in a long paper of Adams and Herb \cite{AH} which establishes such local character identities in a very general setting of coverings of real groups. 


\subsection{ Real groups.}
 In Harish-Chandra's work on the invariant harmonic analysis of real Lie groups, he did not in fact limit himself to the case of linear algebraic groups, but allowed finite central covers of these. Hence, Harish-Chandra's classification of discrete series representations and his Plancherel theorem giving the decomposition of the regular representation $L^2(G)$ were shown for covers of real Lie groups.  Likewise the technique of cohomological induction (the theory of Zuckerman functors) was also developed in this same setting. 
 
 Hence, one understands a lot more about the genuine representation theory of real covering groups. Indeed, there is a classification of such genuine representations called Vogan duality and some representative works in this direction are those of Renard-Trapa \cite{RT1, RT2} and Adams-Trapa \cite{AT}, which led to a Kazhdan-Lusztig algorithm relating the irreducible characters of covering groups and those of standard modules. The recent paper \cite{ABPTV} relates the unitary duals of covering groups and those of an appropriate linear group. 


\subsection{ Invariant harmonic analysis, Eisenstein series and trace formula.}
Our discussion above gives the impression that many results in the representation theory or the theory of automorphic forms on covering groups are based on the study of examples. While this is true to some extent, we would now like to highlight some general results which are necessary ingredients for a systematic theory. 

We begin with invariant harmonic analysis as developed by Harish-Chandra. As mentioned above,  Harish-Chandra's work on the invariant harmonic analysis of real Lie groups applies to finite covers of linear groups, such as his classification of discrete series representations and his Plancherel theorem giving the decomposition of the regular representation $L^2(G)$. His analogous results for $p$-adic groups were written up by Silberger \cite{Si} and also Waldspurger \cite{Wa4}, but only  in the context of linear reductive groups. Recently, many of these foundational results are extended to the covering case (with largely the same proof) by W.-W. Li \cite{L4} (character theory, orbital integrals and Plancherel theorem). Some other results, such as the theory of $R$-groups and the Howe's finiteness conjecture for invariant distributions, were extended to covering groups by C.H. Luo in his thesis work \cite{Lu1}. For smooth representation theory, many standard results developed in \cite{Ca}, such as the Langlands classification and the Casselman square-integrability or temperedness criterion, were extended to general covering groups by Ban-Jantzen \cite{BJ1, BJ2}. It is worth noting that the theory of Bernstein center also works in the same way as in the linear case, as noted by Deligne in his rendition \cite{Be} of the theory of Bernstein center. 


Likewise, in the global setting, the Langlands theory of Eisenstein series was already developed in the setting of covering groups in Moeglin-Waldspurger's monograph \cite{MW}. 
  In a striking series of recent papers \cite{L2, L4, L5, L6}, W.-W. Li has developed the theory of the Arthur-Selberg trace formula for general covering groups, bringing it to the stage of the invariant trace formula.


\section{ A Langlands program for Brylinski-Deligne extensions}
After the historical account of the previous sections, it is natural to ask if the framework of the Langlands program can be extended to the setting of covering groups. The classical Langlands program is built upon the rich and functorial structure theory of linear reductive groups. Such a structure theory of covering groups has now been developed in the work of Brylinski-Deligne and it is our belief that the Brylinski-Deligne theory serves as a good starting point for a systematic extension of the Langlands program to covering groups. 

 \subsection{ What constitutes a Langlands program?}
Before one begins, it may be good to ask what exactly constitutes a Langlands program. For this, one can do no better than to turn to the starting point of the classical Langlands program, which is contained in the famous letter of Langlands to Weil.  The key new ideas introduced in this letter are the notions of the dual group $G^{\vee}$ and the L-group ${^L}G$ of a connected reductive group $\alg{G}$.  Langlands subsequently  reworked in his monograph ``Euler Products" \cite{La1}  the theory of spherical functions and the Satake transform, reinterpreting Satake's results  in the framework of the L-group. This allows him to classify the unramified representations of a quasi-split $p$-adic group in terms of unramified local Galois representations valued in the L-group, which immediately suggests (at least with hindsight) the local Langlands correspondence: classifying all irreducible representations of $\alg{G}(k)$ by local Galois representations valued in ${^L}G$.  This unramified local Langlands correspondence also allows him to introduce the notion of ``automorphic L-functions attached to a finite-dimensional representation of the L-group". 

Thus, a key  ingredient for a Langlands program is undoubtedly the notion of a dual group and an L-group, and a first test for any such candidate dual group or L-group is whether it gives a natural formulation of the Satake isomorphism, leading to a classification of unramified representations. 

A second key realization of Langlands in the initial stage of the classical Langlands program is the difference between conjugacy and stable (or geometric) conjugacy in a reductive group \cite{La2} . More precisely, for a connected reductive group $\alg{G}$ defined over a local field $k$, say, one may consider a coarser equivalence relation than the usual notion of conjugacy in $\alg{G}(k)$. This coarser equivalence relation is conjugacy by elements of $\alg{G}(k^{\sep})$, where $k^{\sep}$ is a separable closure of $k$. This led him to develop the theory of endoscopy, including the definition of endoscopic groups \cite{LL, La3} and  the definition of transfer factors with Shelstad  \cite{LS}.

To summarise, the two key ingredients for a Langlands program are, in our views:
  
  \begin{itemize}
  \item a definition of dual groups and L-groups;
  
  \item a theory of stable conjugacy and endoscopy.
  \end{itemize}
  
  We note that both these ingredients in the classical Langlands program require one to start with a reductive group $\alg{G}$ over $k$, and not just the topological group $\alg{G}(k)$. 
  For example, suppose that $k'/k$ is a separable finite extension of local fields, $\alg{G}$ a reductive group over $k'$ and $\alg{H} : = {\rm Res}_{k'/k} \alg{G}$, so that $\alg{G}(k') = \alg{H}(k)$ as topological groups, and there is no difference between the representation theories of $\alg{G}(k')$ and $\alg{H}(k)$. The dual groups and L-groups of $\alg{G}$ and $\alg{H}$ are however different, even if they can both be used to classify the irreducible representations of the same group. Similarly, the notion of stable conjugacy only makes sense because one has the notion of $k^{\sep}$-points of a reductive group $\alg{G}$ over $k$, with an inclusion $\alg{G}(k) \hookrightarrow \alg{G}(k^{\sep})$. 
This suggests that to have these two ingredients in the setting of  covering groups, one might need to work with covering groups of algebraic origin, such as those provided by the Brylinski-Deligne theory. 
 
 
 \subsection{ Dual Groups}
 We now discuss some prior work on the two key ingredients of a Langlands program highlighted above. As we mentioned in the previous section, people knew what the dual groups of some examples of covering groups should be, such as for $\upr{Mp}_{2n}$, some Kazhdan-Patterson covers and also covers of simply-connected groups. A systematic and general theory was developed in the work of Finkelberg-Lysenko \cite{FL} and Reich \cite{Re} in the context of the Geometric Langlands Program. This was followed in the classical context by the work of McNamara \cite{Mc2} and independently Weissman \cite{W3} who defined the modified dual root datum associated to a Brylinski-Deligne cover, using the invariants associated to such a multiplicative $\alg{K}_2$-torsor by \cite{BD}.

 
 \subsection{ Endoscopy}
The theory of endoscopy for covering groups was initiated by the work of Adams \cite{A1} and Renard \cite{R1, R2} in the context of $\upr{Mp}_{2n}(\RR)$.
The thesis work of J. Schultz considered the case of $\upr{Mp}_2$ over $p$-adic fields.
 The general case of $\upr{Mp}_{2n}$ over any local field was completed in the thesis work of W.-W. Li \cite{L1}, with the endoscopic groups of $\upr{Mp}_{2n}$ being the groups $\upr{SO}_{2a+1} \times \upr{SO}_{2b+1}$, as $(a,b)$ vary over ordered pairs of non-negative integers such that $a+b =n$.   In particular, Li established the transfer of orbital integrals from $\upr{Mp}_{2n}$ to its endoscopic groups,   the fundamental lemma for the unit element of the spherical Hecke algebra and the weighted fundamental lemma \cite{L3}. In his thesis work, C.H. Luo has shown the fundamental lemma for the whole spherical Hecke algebra, as well as established the expected local character identities for the local L-packets of $\upr{Mp}_{2n}$ defined by the local Shimura correspondence of \cite{GS}.
Based on his theory of endoscopy, Li has begun the stabilisation of the invariant trace formula for $\upr{Mp}_{2n}$. It remains to see whether the case of $\upr{Mp}_{2n}$ is an anomaly or is an example of a theory of endoscopy which encompasses a large class of covering groups, such as the Brylinski-Deligne covers. 


 \subsection{ This volume.}
 This brings us to the current volume. 
 
 One of us (M.H.W.) has been thinking about using the Brylinski-Deligne theory as a starting point for the Langlands program for covering groups for some time.  The paper \cite{W1} is an initial attempt to bring the Brylinski-Deligne structure theory to bear on the genuine representation theory of covering tori, whereas the paper \cite{W2} describes the interaction of the Brylinski-Deligne structure theory with the Bruhat-Tits theory of open compact subgroups, answering a question raised at the end of \cite{BD}, while the paper \cite{HW} applies this to the depth zero genuine representation theory of Brylinski-Deligne covers. The paper \cite{W3} gives a definition of the dual group and L-group of a Brylinski-Deligne cover of a split group, using the language of Hopf algebras. This turns out to be overly complicated, making the theory hard to use. Moreover, with hindsight, the candidate L-group there is not always the right one, as it does not make use of all the Brylinski-Deligne invariants.  These initial attempts and ideas were communicated in  a series of letters between M.H.W. and Deligne over the period 2007-2014 and Deligne's ideas and comments have been extremely helpful  in every stage of the development. The interested reader can find this series of letters, which documents the evolution of some of the ideas discussed in this volume, in \cite{W4}.
 
 
  These efforts culminate in the first paper of this volume (by M.H.W.) which defines the L-group of a Brylinski-Deligne cover of a quasi-split group using the language of \'etale gerbes, and tests this L-group for the purpose of representation theory, including the Satake isomorphism and classification of unramified representations as well as the classification of discrete series for covers of real groups. The second paper (by W.T.G. and F.G.) specialises to the case of covers of split groups and introduces another construction of the L-group also due to M.H.W, which is more down-to-earth, as it avoids the language of \'etale gerbes. There is some overlap between the second paper and the first, as the second paper also conducts the necessary tests for the legitimacy of the L-group, namely the Satake isomorphism and the representation theory of covers of split tori. This second paper then goes on to explore some cases of Langlands functoriality such as base change. As a consequence of these two papers, one can now define partial automorphic L-functions for automorphic representations of a Brylinski-Deligne cover. In a followup \cite{Ga} to this work, one of us (F.G.) has extended the results of Langlands' ``Euler Products"  \cite{La1} to the covering setting, using the constant terms of Eisenstein series to show the meromorphic continuation of some of these automorphic L-functions (those of Langlands-Shahidi type).
  Finally, the third paper  (by M.H.W.) of this volume shows that the two notions of L-groups used in the first two papers are in fact the same (for covers of split groups). 
  Since the papers in the volume come with their own extended introductions, we shall refrain from giving a more detailed introduction here.

 Finally, we note that this volume is simply a beginning, and we have only discussed one of the two key ingredients of a Langlands program. We have not addressed the issue of stable conjugacy and endoscopy, except for a brief speculative section in the second paper. We hope that this volume will stimulate further research in this area, leading one day to a fulfilment of the hope expressed by Deligne in his letter  \cite{W4} to M.H.W. (Dated Dec. 14, 2007): 
 
 {\em ``For me, the aim is to understand ``metaplectic" forms on semi-simple groups, the hope
being that they are not ``new" object, but rather correspond to usual automorphic forms
on some other groups, on which they give new information. I would like to have precise
conjectures on the hoped for correspondence, and I view my paper with Brylinski as
setting a landscape in which conjectures should fit. "}

\end{document}

%% file: MathIncludes.tex
\DeclareMathOperator{\K}{K} 
\DeclareMathOperator{\St}{{St}} 
\DeclareMathOperator{\abst}{abs} 
\DeclareMathOperator{\cont}{cont} 
\DeclareMathOperator{\univ}{univ} 

\DeclareMathOperator{\Zar}{Zar}

\DeclareMathOperator{\Hom}{Hom}

\newcommand{\CExt}{CExt}

\DeclareMathOperator{\Ker}{Ker}

\DeclareMathOperator{\Id}{Id}

\DeclareMathOperator{\sep}{sep} 

\newcommand{\From}{\colon}
\newcommand{\Into}{\hookrightarrow}
\newcommand{\Onto}{\twoheadrightarrow}
\newcommand{\To}{\rightarrow}
\newcommand{\isom}{\cong} 

\newcommand{\Cat}[1]{ {\mathsf{#1}} }

\newcommand{\alg}[1]{\boldsymbol{\mathrm{#1}}}
\newcommand{\upr}[1]{ {\mathrm{#1}} }

\newcommand{\ZZ}{\mathbb Z}

\newcommand{\RR}{\mathbb R}
\newcommand{\CC}{\mathbb C}

\renewcommand{\AA}{\mathbb A} 

\renewcommand{\th}{\text{th}}

\newcommand{\Baer}{\dotplus} 

%% file: Asterisque-intro-2.bbl
\begin{thebibliography}{9999999}


\bibitem[A1]{A1}  
J. Adams, {\em Lifting of characters on orthogonal and metaplectic groups}, Duke Math. J. \textbf{92} (1998), no. 1, 129-178.

\bibitem[A2]{A2}  J. Adams, {\em Characters of covering groups of $\upr{SL}(n)$}, J. Inst. Math. Jussieu \textbf{2} (2003), no. 1, 1-21. 

\bibitem[AB]{AB} J. Adams and D. Barbasch, {\em Genuine representations of the
metaplectic group}, Compositio Math. 113 (1998), no. 1, 23-66.
 

\bibitem[ABPTV]{ABPTV}    J. Adams, D. Barbasch, A. Paul, P. Trapa and David Vogan, {\em Shimura correspondences for split real groups}, J. Amer. Math. Soc. \textbf{20} (2007), 701-751.

\bibitem[AH]{AH} J. Adams and R. Herb, {\em Lifting of characters for nonlinear simply laced groups}, 
Represent. Theory \textbf{14} (2010), 70-147.

\bibitem[AT]{AT} J. Adams and P. Trapa, {\em Duality for nonlinear simply laced groups},  Compos. Math. \textbf{148} (2012), no. 3, 931-965.


\bibitem[Ar]{Ar} J. Arthur, \emph{The endoscopic classification of representations: orthogonal and symplectic groups}, 
Colloquium Publications \textbf{61}, American Mathematical Society, 2013.

 

\bibitem[BJ1]{BJ1} D. Ban and C. Jantzen,
\emph{The Langlands quotient theorem for finite central extensions of p-adic groups},
Glas. Mat. Ser. III \textbf{48} (68) (2013), no. 2, 313-334.

\bibitem[BJ2]{BJ2} D. Ban and C. Jantzen,
{\em The Langlands quotient theorem for finite central extensions of p-adic groups II: intertwining operators and duality},  Glas. Mat. Ser. III \textbf{51} (71) (2016), no. 1, 153-163.

\bibitem[BLS]{BLS} W. D. Banks, J. Levy and M. R.  Sepanski, {\em  Block-compatible metaplectic cocycles},  J. Reine Angew. Math. \textbf{507} (1999), 131-163.

 
\bibitem[Be]{Be} J. Bernstein, {\em Le ``centre" de Bernstein}, edited by P. Deligne, in Travaux en Cours, Representations of reductive groups over a local field, 1-32, Hermann, Paris, 1984.

 
 
 \bibitem[BBL]{BBL} W. Banks, D.  Bump, D. Lieman, {\em  Whittaker-Fourier coefficients of metaplectic Eisenstein series},  Compositio Math. \textbf{135} (2003), no. 2, 153-178. 

   \bibitem[BBCFG]{BBCFG}
B. Brubaker, D. Bump, G. Chinta, S. Friedberg and P. Gunnells,
\emph{Metaplectic ice} in \emph{Multiple Dirichlet series, L-functions and automorphic forms}
Progr. Math., 300, Birkhauser, 2012, 65--92. 



\bibitem[BBF]{BBF}
B. Brubaker, D. Bump and S. Friedberg,
\emph{Weyl group multiple Dirichlet series, Eisenstein series and crystal bases},
Ann. of Math. \textbf{173} (2011), no. 2, 1081-1120.

\bibitem[BBFH]{BBFH}
B. Brubaker, D. Bump, S. Friedberg and J. Hoffstein,
\emph{Weyl group multiple Dirichlet series. III. Eisenstein series and twisted unstable $A_r$},
Ann. of Math. \textbf{166} (2007), No. 1, 293-316.

\bibitem[BrF]{BrF} B.  Brubaker and S.  Friedberg,
\emph{Whittaker coefficients of metaplectic Eisenstein series},
Geom. Funct. Anal. Vol. \textbf{25} (2015), 1180-1239.


\bibitem[BrFH]{BrFH} B.  Brubaker, S.  Friedberg and J.  Hoffstein, {\em Cubic twists of $\upr{GL}(2)$ automorphic L-functions}, Invent. Math. \textbf{160} (2005), no. 1, 31-58.

\bibitem[BD]{BD}
J.-L. Brylinski and P. Deligne,
\emph{Central extensions of reductive groups by $\K_2$},
Publ. Math. Inst. Hautes $\acute{\text{E}}$tudes Sci. \textbf{94} (2001), 5-85.

  


\bibitem[BFH1]{BFH1} D. Bump, S. Friedberg, J. Hoffstein, {\em Eisenstein series on the metaplectic group and nonvanishing theorems for automorphic L-functions and their derivatives},  Ann. of Math. (2) \textbf{131} (1990), no. 1, 53-127.

\bibitem[BFH2]{BFH2} D. Bump, S. Friedberg, J. Hoffstein, {\em  Nonvanishing theorems for L-functions of modular forms and their derivatives},  Invent. Math. \textbf{102} (1990), no. 3, 543-618. 
 
 \bibitem[BFH3]{BFH3} D. Bump, S. Friedberg, J.  Hoffstein, {\em  p-adic Whittaker functions on the metaplectic group},  Duke Math. J. \textbf{63} (1991), no. 2, 379-397.
 
\bibitem[BG]{BG} D. Bump and D. Ginzburg, {\em  Symmetric square L-functions on $\upr{GL}(r)$},  Ann. of Math. (2) \textbf{136} (1992), no. 1, 137-205. 
 
 
\bibitem[BF]{BF} D. Bump and S. Friedberg, {\em Metaplectic generating functions and Shimura integrals},  Automorphic forms, automorphic representations, and arithmetic (Fort Worth, TX, 1996), 1-17, Proc. Sympos. Pure Math., 66, Part 2, Amer. Math. Soc., Providence, RI, 1999. 
 

\bibitem[BH1]{BH1} D. Bump and J. Hoffstein, {\em Cubic metaplectic forms on $\upr{GL}(3)$}, Invent. Math. \textbf{84} (1986), no. 3, 481-505.

\bibitem[BH2]{BH2} D. Bump, J.  Hoffstein, {\em Some Euler products associated with cubic metaplectic forms on $\upr{GL}(3)$}, Duke Math. J. \textbf{53} (1986), no. 4, 1047-1072. 
 

\bibitem[BH3]{BH3} D. Bump, J.  Hoffstein,  {\em  On Shimura's correspondence},  Duke Math. J. 55 (1987), no. 3, 661-691. 
  
\bibitem[BH4]{BH4} D. Bump, J.  Hoffstein,  {\em   Some conjectured relationships between theta functions and Eisenstein series on the metaplectic group},  Number theory (New York, 1985/1988), 1-11, Lecture Notes in Math., 1383, Springer, Berlin, 1989.

 
 
 \bibitem[C1]{C1} Y.Q. Cai, {\em Fourier Coefficients for Theta Representations on Covers of General Linear Groups}, available  at \url{https://arxiv.org/abs/1602.06614}.
 
\bibitem[C2]{C2}   Y.Q. Cai, {\em Fourier Coefficients for Degenerate Eisenstein Series and the Descending Decomposition}, available  at \url{https://arxiv.org/abs/1606.08872}.
 
 \bibitem[CFGK]{CFGK} Y. Q. Cai, S. Friedberg, D. Ginzburg and E. Kaplan, {\em Doubling Constructions for Covering Groups and Tensor Product L-Functions}, available  at \url{http://arxiv.org/abs/1601.08240}.
 
  \bibitem[Ca]{Ca} W. Casselman, {\em Introduction to the theory of admissible representations for reductive $p$-adic groups,} Draft 7 May 1993.
  
 \bibitem[CG1]{CG1} G. Chinta and P. Gunnells, {\em Weyl group multiple Dirichlet series constructed from quadratic characters},  Invent. Math. \textbf{167} (2007), no. 2, 327-353.
  
  \bibitem[CG2]{CG2} G. Chinta and P. Gunnells,  {\em  Constructing Weyl group multiple Dirichlet series},  J. Amer. Math. Soc. \textbf{23} (2010), no. 1, 189-215.
  
  
\bibitem[CO]{CO}
G. Chinta and O. Offen,
\emph{A metaplectic Casselman-Shalika formula for $\upr{GL}_r$},
Amer. J. Math., \textbf{135} (2013), 403-441.


\bibitem[D1]{D1} P. Deligne, {\em Sommes de Gauss cubiques et rev\^{e}tements de $\upr{SL}(2)$ (d'apr\`{e}s S. J. Patterson)},  S\'eminaire Bourbaki (1978/79), Exp. No. 539, pp. 244-277, 
Lecture Notes in Math., 770, Springer, Berlin, 1980.


\bibitem[D2]{D2} P. Deligne, {\em Extensions centrales de groupes alg\'ebriques simplement connexes et cohomologie galoisienne}, 
Publ. Math. Inst. Hautes $\acute{\text{E}}$tudes Sci. (1996), no. 84, 35-89.

\bibitem[D3]{D3} P. Deligne, {\em Expos\'e XVIII:  La formule de dualit\'e globale}, SGA IV, tome 3.

\bibitem[De]{De} V. Deodhar, {\em On central extensions of rational points of algebraic groups}, 
 Amer. J. Math. \textbf{100} (1978), no. 2, 303-386. 

  \bibitem[FL]{FL}  M. Finkelberg and S. Lysenlo, {\em Twisted geometric Satake equivalence}, 
J. Inst. Math. Jussieu \textbf{9} (2010), no. 4, 719-739.

\bibitem[F]{F} Y. Flicker, {\em Automorphic forms on covering groups of $\upr{GL}(2)$}, Invent. Math. \textbf{57} (1980), no. 2, 119-182.


\bibitem[FG1]{FG1} S. Friedberg and D. Ginzburg,
{\em On the genericity of Eisenstein series and their residues for covers of $GL_m$}, to appear in Int. Math. Res. Not..

\bibitem[FG2]{FG2} S. Friedberg and D. Ginzburg, 
{\em Descent and theta functions for metaplectic groups}, 
to appear in the J. Euro. Math. Soc..

\bibitem[FG3]{FG3} S. Friedberg and D. Ginzburg,
{\em Theta functions on covers of symplectic groups}, to appear in {\em Automorphic Forms, L-functions and Number Theory}, a special issue of BIMS in honor of Prof. Freydoon ShahidiÕs 70th birthday.


\bibitem[FH]{FH} S. Friedberg and J. Hoffstein, {\em Nonvanishing theorems for automorphic L-functions on $\upr{GL}(2)$}. Ann. of Math. (2) \textbf{142} (1995), no. 2, 385-423.

\bibitem[FZ]{FZ} S. Friedberg and L. Zhang, {\em Eisenstein series on covers of odd orthogonal groups}, Amer. J. Math. \textbf{137} (2015), 953Ð1011.


\bibitem[FK]{FK}
Y. Flicker and D. Kazhdan,
\emph{Metaplectic correspondence},
Publ. Math. Inst. Hautes $\acute{\text{E}}$tudes Sci. \textbf{64} (1986), 53-110.

\bibitem[G1]{G1}
W.~T.~Gan,
\emph{The Shimura correspondence \`a la Waldspurger},
available at \url{http://www.math.nus.edu.sg/~matgwt/postech.pdf}.

\bibitem[G2]{G2} W.T. Gan, {\em Theta correspondence: recent progress and applications}, Procedings of ICM (2014) Vol. 2, 343-366.

\bibitem[GGP]{GGP}  W.T. Gan, B.H. Gross and D.Prasad, {\em Symplectic local root numbers, central critical L-values and restriction problems in the representation theory of classical groups}, Asterisque 346, 1-110.

\bibitem[GS]{GS}
W.~T.~Gan and G.~Savin,
\emph{Representations of metaplectic groups I:
 epsilon dichotomy and local Langlands correspondence},
Compositio Math. \textbf{148} (2012), 1655--1694.

\bibitem[Ga]{Ga} F. Gao,
{\em The Langlands-Shahidi L-functions for Brylinski-Deligne extensions}, to appear in Amer. J. Math.. 
 

\bibitem[HP]{HP} D. R. Heath-Brown and S.J.  Patterson, {\em  The distribution of Kummer sums at prime arguments},  J. Reine Angew. Math. \textbf{310} (1979), 111-130.
 
\bibitem[HW]{HW} T. K. Howard and M. H, Weissman, {\em Depth-zero representations of nonlinear covers of p-adic groups}, Int. Math. Res. Not. (2009), no. 21, 3979-3995.


\bibitem[KP1]{KP1}
D. A. Kazhdan and S. J. Patterson,
\emph{Metaplectic forms},
Publ. Math. Inst. Hautes $\acute{\text{E}}$tudes Sci. \textbf{59} (1984), 35-142.

\bibitem[KP2]{KP2}
D. A. Kazhdan and S. J. Patterson,
\emph{Towards a generalized Shimura correspondence},
Adv. in Math. \textbf{60(2)} (1986), 161-234.



\bibitem[Ka]{Ka} A. Kable, {\em The tensor product of exceptional representations on the general linear group},  Ann. Sci. cole Norm. Sup. (4) \textbf{34} (2001), no. 5, 741-769. 

\bibitem[Ku1]{Ku1} T. Kubota, {\em Topological covering of $\upr{SL}(2)$ over a local field},  J. Math. Soc. Japan \textbf{19} (1967), 114-121.


\bibitem[Ku2]{Ku2} T. Kubota, {\em On automorphic functions and the reciprocity law in a number field}, 
Lectures in Mathematics, Department of Mathematics, Kyoto University, No. 2 Kinokuniya Book-Store Co., Ltd., Tokyo 1969 iii+65 pp. 

\bibitem[LL]{LL} J.P. Labesse and R. Langlands, {\em  L-indistinguishability for $\upr{SL}(2)$},  Canad. J. Math. \textbf{31} (1979), no. 4, 726-785.

\bibitem[La1]{La1}
 R. Langlands, \emph{Euler products},
 Yale Univ. Press, 1971.

\bibitem[La2]{La2}
 R. Langlands, {\em  Stable conjugacy: definitions and lemmas},  Canad. J. Math. \textbf{31} (1979), no. 4, 700-725. 
 
\bibitem[La3]{La3}
 R. Langlands,  {\em Les d\'ebuts d'une formule des traces stable}, Publications Math\'ematiques de l'Universit\'e Paris VII,  U.E.R. de MathÂmatiques, Paris, 1983. v+188 pp.

 

\bibitem[LS]{LS} R. Langlands and D. Shelstad, {\em On the definition of transfer factors},  Math. Ann. \textbf{278} (1987), no. 1-4, 219-271.


\bibitem[LM1]{LM1} E. Lapid and Z.Y. Mao, {\em A conjecture on Whittaker-Fourier coefficients of cusp forms}, J. Number Theory \textbf{146} (2015), 448-505.

\bibitem[LM2]{LM2} E. Lapid and Z.Y. Mao, {\em Whittaker-Fourier coefficients of cusp forms on $\widetilde{\upr{Sp}}_n$: reduction to a local statement}, available  at \url{https://arxiv.org/abs/1401.0198}.

 \bibitem[LM3]{LM3} E. Lapid and Z.Y. Mao, {\em  On an analogue of the Ichino--Ikeda conjecture for Whittaker coefficients on the metaplectic group}, available  at \url{https://arxiv.org/abs/1404.2905}.


\bibitem[L1]{L1} W.-W. Li, {\em Transfert d'int\'egrales orbitales pour le groupe m\'etaplectique}, Compositio Math. \textbf{147} (2011), 524-590.

\bibitem[L2]{L2} W.-W. Li, {\em La formule des traces pour les rev\^etements de groupes r\'eductifs connexes. I. Le d\'eveloppement 
g\'eom\'etrique fin}, J. Reine Angew. Math. \textbf{686} (2014), 37-109. 

\bibitem[L3]{L3} W.-W. Li, {\em Le lemme fondamental pond\'er\'e pour le groupe m\'etaplectique}, 
Canadian J. Math. Vol \textbf{64 (3)} (2012), 497-543.

\bibitem[L4]{L4} W.-W. Li, {\em La formule des traces pour les rev\^etements de groupes r\'eductifs connexes. II. 
Analyse harmonique locale}, Annales scientifiques de l'ENS \textbf{45} (2012), 787-859.

\bibitem[L5]{L5} W.-W. Li, {\em La formule des traces pour les rev\^etements de groupes r\'eductifs connexes. III. 
Le d\'eveloppement spectral fin}, Math. Ann. \textbf{356 (3)} (2013), 1029-1064.

\bibitem[L6]{L6} W.-W. Li, {\em La formule des traces pour les rev\^etements de groupes r\'eductifs connexes. IV. 
Distributions invariantes},
Ann. de l'Institut Fourier \textbf{64} No. 6 (2014), 2379-2448.

\bibitem[L7]{L7} 
W.-W. Li, 
{\em La formule des traces stable pour le groupe m\'etaplectique: les termes elliptiques},
Invent. Math. \textbf{202} (2015), 743-838.

 \bibitem[Lu1]{Lu1} C.H. Luo, {\em  Howe conjecture for covering groups}, preprint (2017).
 
 \bibitem[Lu2]{Lu2} C.H. Luo, {\em   Spherical fundamental lemma for metaplectic groups}, preprint (2017).
 
 \bibitem[Lu3]{Lu3} C.H. Luo, {\em   Endoscopic character identities for metaplectic groups}, preprint (2017)
 
\bibitem[Ma]{Ma}
H. Matsumoto,
\emph{Sur les sous-groupes arithm\'etiques des groupes semi-simples d\'eploy\'es},
Ann. Sci. \'Ecole Norm. Sup. \textbf{4} (1969), 1-62.

\bibitem[Mac]{Mac}
George W. Mackey, {\em Les ensembles Bor\'elien et les extensions des groupes}, J. Math. Pures Appl. \textbf{36} (1957), 171-178.


\bibitem[Mc1]{Mc1}
P. McNamara, 
{\em Metaplectic Whittaker functions and crystal bases},  Duke Math. J. \textbf{156} (2011), no. 1, 1-31.

\bibitem[Mc2]{Mc2}
P. McNamara,
\emph{Principal series representations of metaplectic groups over local fields}
in \emph{Multiple Dirichlet series, $L$-functions and automorphic forms},
Birkhauser, 2012, 299-328.

\bibitem[Mc3]{Mc3}
P. McNamara, {\em The Metaplectic Casselman-Shalika Formula},
Trans. Amer. Math. Soc. \textbf{368} No. 4 (2016), 2913-2937.

\bibitem[Me1]{Me1} P. Mezo, {\em Comparisons of general linear groups and their metaplectic coverings I}, Canad. J. Math. \textbf{54} (2002), no. 1, 92-137.

\bibitem[Me2]{Me2} P. Mezo, {\em Comparisons of general linear groups and their metaplectic coverings II},  Represent. Theory \textbf{5} (2001), 524-580.

\bibitem[Mi]{Mi} J. Milnor, {\em Introduction to algebraic K-theory}, Annals of Mathematics Studies, No. 72. Princeton University Press, Princeton, N.J.; University of Tokyo Press, Tokyo, 1971. xiii+184 pp

\bibitem[Mo]{Mo}
C. Moore,
\emph{Group extensions of $p$-adic and adelic linear groups},
Publ. Math. Inst. Hautes $\acute{\text{E}}$tudes Sci. \textbf{35} (1968), 157-222.


\bibitem[Mo1-2]{Mo1-2} C. Moore, {\em Extensions and low dimensional cohomology theory of locally compact groups. I, II}, Trans. Amer. Math. Soc. \textbf{113} (1964) 40-63.

\bibitem[Mo3]{Mo3}  C. Moore, {\em Group extensions and cohomology for locally compact groups. III},  Trans. Amer. Math. Soc. \textbf{221} (1976), no. 1, 1-33.

\bibitem[Mo4]{Mo4}  C. Moore, {\em Group extensions and cohomology for locally compact groups. IV}, Trans. Amer. Math. Soc. \textbf{221} (1976), no. 1, 35-58.

\bibitem[MW]{MW}
 C. Moeglin and J.-L. Waldspurger,
 \emph{Spectral decomposition and Eisenstein series},
 Cambridge University Press, 1995.



\bibitem[P1]{P1} S. J. Patterson, {\em A cubic analogue of the theta series}, J. Reine Angew. Math. \textbf{296} (1977), 125-161.

\bibitem[P2]{P2} S. J. Patterson, {\em A cubic analogue of the theta series II},  J. Reine Angew. Math. \textbf{296} (1977), 217-220.

\bibitem[Pr]{Pr} D. Prasad, {\em Notes on central extensions}, available at \url{https://arxiv.org/abs/1502.02140}.

\bibitem[P]{P} G. Prasad, {\em Deligne's topological central extension is universal}, Adv. Math. 181 (2004), no. 1, 160-164.

\bibitem[PR1]{PR1} G. Prasad and M.S. Raghunathan, {\em Topological central extensions of semisimple groups over local fields},
 Ann. of Math. (2) \textbf{119} (1984), no. 1, 143-201.

\bibitem[PR2]{PR2} G. Prasad and M.S. Raghunathan, {\em Topological central extensions of semisimple groups over local fields. II}, 
 Ann. of Math. (2) \textbf{119} (1984), no. 2, 203-268. 

\bibitem[PR3]{PR3} G. Prasad and M.S. Raghunathan,
 {\em Topological central extensions of $\upr{SL}_1(D)$},  Invent. Math. \textbf{92} (1988), no. 3, 645-689.


\bibitem[PR]{PR} G. Prasad and A. Rapinchuk, {\em Computation of the metaplectic kernel},  Inst. Hautes \'Etudes Sci. Publ. Math. No. \textbf{84} (1996), 91-187.

\bibitem[Ra]{Ra}
R. Rao,
\emph{On some explicit formulas in the theory of Weil representation},
Pacific J. Math. \textbf{157} (1993), No. 2, 335-371.

\bibitem[Re]{Re} R.C. Reich,
\emph{Twisted geometric Satake equivalence via gerbes on the factorizable Grassmannian},
Represent. Theory \textbf{16} (2012), 345-449.

\bibitem[R1]{R1} D. Renard, {\em Transfert d'int\'{e}grales orbitales entre $\upr{Mp}(2n,\RR)$ et $\upr{SO}(n+1,n)$},  Duke Math. J. \textbf{95} (1998), no. 2, 425-450.

\bibitem[R2]{R2} D. Renard, {\em Endoscopy for $\upr{Mp}(2n,\RR)$},  Amer. J. Math. \textbf{121} (1999), no. 6, 1215-1243.

\bibitem[RT1]{RT1} D. Renard and P. Trapa, {\em Irreducible genuine characters of the metaplectic group: Kazhdan-Lusztig algorithm and Vogan duality},  Represent. Theory \textbf{4} (2000), 245-295. 

\bibitem[RT2]{RT2} D. Renard and P. Trapa, {\em Kazhdan-Lusztig algorithms for nonlinear groups and applications to Kazhdan-Patterson lifting}, Amer. J. Math. \textbf{127} (2005), no. 5, 911-971. 

 \bibitem[Sa1]{Sa1}
G. Savin,  {\em Local Shimura correspondence}, Math. Ann. \textbf{280} (1988), no. 2, 185-190. 
 
\bibitem[Sa2]{Sa2}
G. Savin,
\emph{On unramified representations of covering groups},
J. Reine Angew. Math. \textbf{566} (2004), 111-134.

\bibitem[Sc]{Sc} 
J. Schultz, {\em Lifting of characters of $\tilde{\upr{SL}}_2(F)$ and $\upr{SO}_{1,2}(F)$ for $F$ a nonarchimedean local field}, PhD thesis, Univ. of Maryland, 1998.

\bibitem[S1]{S1}
R. Steinberg,
\emph{G\'en\'erateurs, relations et rev\^etements de groupes alg\'ebriques}, in
\emph{Colloq. Th\'eorie des Groupes Alg\'ebriques} (Bruxelles, 1962), Louvain, 1962.

\bibitem[S2]{S2}
R. Steinberg,
\emph{Lectures on Chevalley groups}, 
notes prepared by John Faulkner and Robert Wilson, Yale University, New Haven, Conn., 1968.


\bibitem[Sh]{Sh}
G.~Shimura,
\emph{On modular forms of half integral weight}, 
Ann. of Math. (2) \textbf{97} (1973), 440--481. 

\bibitem[Si]{Si}
A. Silberger,
\emph{Introduction to harmonic analysis on reductive p-adic groups}. Based on lectures by Harish-Chandra at the Institute for Advanced Study, 1971-1973. Princeton University Press, 1979.

\bibitem[Su1]{Su1}
T. Suzuki,
\emph{Rankin-Selberg convolution of generalized theta series},
J. reine angew. Math. \textbf{414} (1991), 149-205.

\bibitem[Su2]{Su2}
T. Suzuki,
\emph{Metaplectic Eisenstein series and the Bump-Hoffstein conjecture},
Duke Math. J. \textbf{90} No. 3 (1997), 577-630.



\bibitem[Sz1]{Sz1}
D. Szpruch,
\emph{The Langlands-Shahidi method for the metaplectic group and applications},
thesis (Tel Aviv University), available at \url{arXiv:1004.3516v1}.

\bibitem[Sz2]{Sz2} D. Szpruch, {\em  Some irreducibility theorems of parabolic induction on
 the metaplectic group via the Langlands-Shahidi method}, Israel J. Math. \textbf{195} (2013), no. 2, 897-971.





\bibitem[Ta]{Ta} S. Takeda, {\em The twisted symmetric square L-function of $GL(r)$}, Duke Math. J. \textbf{163} (2014), 175-266.


\bibitem[Tu]{Tu} J. Tunnell, {\em A classical Diophantine problem and modular forms of weight $3/2$},  Invent. Math. \textbf{72} (1983), no. 2, 323-334.

\bibitem[Wa1]{Wa1} J.L. Waldspurger, {\em Correspondance de Shimura}, J. Math. Pures Appl. (9) \textbf{59} (1980), no. 1, 1-132.

\bibitem[Wa2]{Wa2} J.L. Waldspurger, {\em  Sur les coefficients de Fourier des formes modulaires de poids demi-entier},  J. Math. Pures Appl. (9) \textbf{60} (1981), no. 4, 375-484.
 
\bibitem[Wa3]{Wa3} J.L. Waldspurger, {\em Correspondances de Shimura et quaternions},  Forum Math. 3 (1991), no. 3, 219-307.
 
\bibitem[Wa4]{Wa4} J.L. Waldspurger, {\em La formule de Plancherel pour les groupes p-adiques (d'apr\`es Harish-Chandra)},  J. Inst. Math. Jussieu 2 (2003), no. 2, 235-333. 




\bibitem[We1]{We1} A. Weil, {\em Sur certains groupes d'op\'erateurs unitaires},  Acta Math. \textbf{111} (1964), 143-211.

\bibitem[We2]{We2} A. Weil, {\em Sur la formule de Siegel dans la th\'eorie des groupes classiques}, Acta Math. \textbf{113} (1965),  1-87. 

\bibitem[W1]{W1} M.~Weissman, \emph{Metaplectic tori over local fields}, Pacific J. Math., \textbf{241 (1)} (2009), 169-200.

\bibitem[W2]{W2} M.~Weissman,
\emph{Managing metaplectiphobia: covering $p$-adic groups} in
\emph{Harmonic analysis on reductive, $p$-adic groups} 2011, 237-277.

\bibitem[W3]{W3}
M.~Weissman,
\emph{Split metaplectic groups and their $L$-groups},
J. Reine Angew. Math. \textbf{696} (2014), 89-141.

\bibitem[W4]{W4} M.~Weissman, {\em Correspondence with P. Deligne}, available at \url{http://martyweissman.com/WeissmanDeligneLetters.pdf}.


 







\end{thebibliography}
